\documentclass[a4paper,12pt,reqno]{article}
\pdfoutput=1 %To force HAL and Arxiv to use pdflatex

\usepackage{graphicx}
\usepackage{tabularx}
\usepackage{amsfonts}
\usepackage{amssymb}
\usepackage{pstricks}
\usepackage{amsmath,upgreek}
\usepackage{array}
\usepackage{pgfplots}
\usepackage{pstricks}
\usepackage{bm}
\usepackage[top=3cm, bottom=3cm, left=3cm, right=3cm]{geometry}
\pgfplotscreateplotcyclelist{ageplot}{%
{blue,mark=o, mark size={3}},
{red,mark=square},
{black,mark=star},
{brown,mark=triangle},
{green,mark=diamond},}

\tikzset{every mark/.append style={scale=0.7}}
\usepackage{pgfplotstable,booktabs}
\pgfplotstableset
{% global config, for example in the preamblethese columns/<colname>/.style={<options>} things define a style which applies to <colname> only.
empty cells with={--}, % replace empty cells with '--'
every head row/.style={before row=\toprule,after row=\midrule},
every last row/.style={after row=\bottomrule}
}

\newcommand\vm[1]{\bm{\mathrm{#1}}} 		% Vector or matrix
\newcommand{\norm}[1]{\lVert#1\rVert} 		% Norm
\newcommand{\abs}[1]{\left|#1\right|} 		% Absolute value
\newcommand{\dd}{\rm d} 			        % Derivative
\title{Locally equilibrated stress recovery for goal oriented error estimation in the extended finite element method}

\author{O.A. Gonz\'alez-Estrada$^{*1}$ \and J.J. R\'odenas$^{2}$ \and  S.P.A. Bordas$^{1}$ \and E. Nadal$\,^{2}$ \and  P. Kerfriden$^{1}$ \and F.J. Fuenmayor$\,^{2}$}

\begin{document}

\maketitle

\begin{center}\small{
$^{1}$Institute of Mechanics and Advanced Materials (IMAM), Cardiff School of Engineering, Cardiff University, 
Queen's Buildings, The Parade, Cardiff CF24 3AA Wales, UK, \\
$^{2}$Centro de Investigaci\'on de Tecnolog\'ia de Veh\'iculos (CITV), \\ Universitat Polit\`{e}cnica de Val\`{e}ncia, E-46022-Valencia, Spain}
\end{center}

%e-mail: \{estradaoag, kerfridenp\}@cardiff.ac.uk, stephane.bordas@alum.northwestern.edu   , e-mail: jjrodena@mcm.upv.es, ennaso@upv.es, ffuenmay@mcm.upv.es

%\corraddr{J. J. R\'odenas, Centro de Investigaci\'on de Tecnolog\'ia de Veh\'iculos, Departamento de Ingenier\'ia Mec\'anica y de Materiales, Universitat Polit\`ecnica de Val\`encia, Camino de Vera s/n 46022-Valencia (Spain). E-mail: jjrodena@mcm.upv.es}

\begin{abstract}
 Goal oriented error estimation and adaptive procedures are essential for the accurate and efficient evaluation of finite element numerical simulations that involve complex domains. By locally improving the approximation quality, for example, by using the extended finite element method (XFEM), we can solve expensive problems which could result intractable otherwise. Here, we present an error estimation technique for enriched finite element approximations that is based on an equilibrated recovery technique, which considers the stress intensity factor as the quantity of interest. The locally equilibrated superconvergent patch recovery is used to obtain enhanced stress fields for the primal and dual problems defined to evaluate the error estimate.
\end{abstract}

\vspace{10pt}

%TODO check variation with respect to radius

%==============================================================================
\section{Introduction}
%==============================================================================

Nowadays, complex mechanical problems are solved using large numerical simulations in many engineering settings. One particular aspect of the design process is to offer good reliable solutions with the lowest computational cost. As numerical methods introduce an error in the solution due to the approximations used to solve the problem, it becomes necessary to quantify this error in order to guarantee the quality of the results \cite{ainsworthoden2000, ladevezepelle2005}. Moreover, in order to increase the computer efficiency, it is common practice to use adaptivity procedures to improve the accuracy whilst keeping the model with a tractable small size. 

Since the beginning of the use of numerical simulations many methods have been developed to control the discretisation error of finite element approximations, mostly based on the evaluation of global energy norms. These methods can be broadly classified in residual based \cite{babuskarheinboldt1978}, recovery based \cite{zienkiewiczzhu1987} and dual analysis \cite{pereiraalmeida1999}. However, a more interesting approach is to control the error in a particular quantity relevant for the design process \cite{ainsworthoden2000,paraschivoiuperaire1997, ladevezerougeot1999, odenprudhomme2001}. This quantity could be defined as a bounded functional that describes the displacement or stresses in a given area of the domain, or for the case of fracture mechanics, the stress intensity factor that characterises the crack. This approach, referred to as goal oriented, is usually based on the use of duality techniques that involve the formulation of an \emph{adjoint} or \emph{dual} problem directly related to the quantity of interest (QoI). Residual methods have been frequently used to evaluate the error in quantities of interest although examples involving recovery techniques  can be found in \cite{cirakramm1998,ruterstein2006}, and considering dual analysis in \cite{almeidapereira2006}. In \cite{ruterstein2006}, recovery and residual based estimates of the error in evaluating the $J$-integral for finite element (FE) approximations  in the context of linear elastic fracture mechanics are presented. The numerical results showed that a recovery technique with a standard superconvergent patch recovery (SPR) gives more accurate results than the residual estimates presented. Note that if we use an energy estimate with bounding properties, then the error estimate for the quantity of interest is bounded \cite{odenprudhomme2001, ladevezerougeot1999}. On the other hand, it is usually difficult to obtain guaranteed error bounds of the quantities of interest while maintaining the accuracy of the estimate. The need of such a bound is also arguable in an engineering context as the reliability of an a posteriori error estimate, which is quantified by its local effectivity, can be verified beforehand on a number of practical cases. Here, we are interested in increasing the effectivity of the error estimate used to guide adaptive algorithms rather than error bounding.

In the context of fracture mechanics, the extended finite element method (XFEM) \cite{moesdolbow1999} has been successfully used to enrich the finite element approximation in order to represent the particular features of cracks, namely, the discontinuity along the crack faces and the singularity at the crack tip. This method helps to overcome some of the difficulties when modelling crack propagation, such as the need for remeshing to obtain conforming meshes to the crack topology. Error estimators in energy norm for XFEM and other partition of unity methods have been proposed in \cite{ bordasduflot2007, duflotbordas2008, rodenasgonzalez2008, rodenasgonzalez2010} using recovery techniques, and in \cite{stroubouliszhang2006, panetierladeveze2010,gerasimovruter2012} using the residual approach. 
A goal oriented approach for enriched finite element approximations based on the constitutive relation error has been presented in \cite{panetierladeveze2009}. In \cite{rutergerasimov2013} goal oriented error estimators based on the explicit residual method were introduce for the XFEM framework. In \cite{pannachetsluys2009}, adaptive techniques based on energy norm and goal oriented error estimation have been investigated for enriched finite element approximations.

%In this paper, we propose a goal oriented error estimation technique for XFEM approximations that is based in the enhanced recovery technique previously presented in \cite{rodenasgonzalez2008, rodenasgonzalez2010}. We use the stress intensity factor (SIF) typical of fracture mechanics problems as the quantity of interest. As shown in \cite{bordasduflot2007,gonzalezrodenas2012}, the use of standard recovery techniques (e.g. SPR) for XFEM approximations of fracture mechanics problems does not give good results and is not recommended for accurate error estimates in energy norm. The use of enhanced recovery techniques is recommended in these references. Therefore,  error estimates in quantities of interest will also require a careful consideration of the singular character of the XFEM solution, and the use of extended recovery approaches becomes a necessity to obtain accurate estimates. To improve the quality of the recovered stresses for the primal and dual problems we consider equilibrium constraints locally at patches and the splitting of the stress field to describe the singular behaviour of the solution. 

In this paper, we propose a goal oriented error estimation technique for XFEM approximations that is based on the enhanced recovery technique previously presented in \cite{rodenasgonzalez2008, rodenasgonzalez2010}. We use the stress intensity factor (SIF) typical of fracture mechanics problems as the quantity of interest. As shown in \cite{bordasduflot2007, gonzalezrodenas2012}, error estimators based on standard recovery techniques (e.g. SPR) provide inaccurate results because the polynomial basis of the recovered stress field is unable to improve the XFEM solution in fracture mechanics problems, which includes the singular terms. The use of enhanced recovery techniques is recommended in these references. Therefore,  error estimates in quantities of interest will also require a careful consideration of the singular character of the XFEM solution, and the use of extended recovery approaches becomes a necessity to obtain accurate estimates. To improve the quality of the recovered stresses for the primal and dual problems, and therefore, the accuracy of the error estimate, we consider equilibrium constraints locally in patches of elements and the splitting of the stress field to describe the singular behaviour of the solution. 

The paper is organised as follows.  In Section  \ref{sec:ProbStatement}, we introduce the problem under consideration and its corresponding enriched approximation. The general framework for error measures is presented in Section \ref{sec:error}. In Section \ref{sec:errorQoI}, we show useful analytical definitions of QoI for the enforcement of equilibrium conditions. We discuss the formulation of the dual problem when considering the stress intensity factor as the quantity of interest in the goal oriented approach. Numerical results are provided in Section \ref{sec:results} and conclusion are drawn in Section \ref{sec:conclusions}.

%==============================================================================
\section{Problem statement and XFEM solution}
%==============================================================================
\label{sec:ProbStatement}

In this section, we introduce the 2D linear elasticity problem. We denote by $\vm{u}$ the displacement, by $\vm{\sigma}$ the Cauchy stress  and by $\vm{\varepsilon}$ the strain, all these fields defined over the domain $\Omega \subset \mathbb{R}^{2}$, of boundary denoted by $\partial \Omega$.  $\Gamma _{N}$ and $\Gamma _{D}$ refer to the parts of the boundary where the Neumann and Dirichlet conditions are applied, and $\Gamma_C$ to the free traction surface describing a crack such that $\partial \Omega = \Gamma_N \cup \Gamma_D\cup \Gamma_C$ and $\Gamma_N \cap \Gamma_D \cap \Gamma_C =\emptyset$. We denote as $\vm{b}$ the body loads, $\vm{t}$ the tractions imposed along $\Gamma _N$ and $\vm{\sigma}_0$, $\vm{\varepsilon}_0$ the initial stresses and strains. The displacement field $\vm{u}$ is the  solution of the problem given by 
\begin{align}
  \vm{L}^T \vm{\sigma} + \vm{b} &=  0    &&  {\rm in }\; \Omega   \label{Eq:IntEq} \\
  \vm{G} \vm{\sigma}  &= \vm{t}     &&  {\rm on }\; \Gamma _{N}                       \label{Eq:Neumann}\\
  \vm{G} \vm{\sigma}   &= \vm{0}     &&  {\rm on }\; \Gamma _{C}\\
  \vm{u}  &= {\vm{0}} &&  {\rm on }\; \Gamma _{D} 
  \label{Eq:Dirichlet} \\
  \vm{\varepsilon} (\vm{u}) &= \vm{L} \vm{u}  &&  {\rm in }\; \Omega \\
  \vm{\sigma} &= \vm{D} ( \vm{\varepsilon}(\vm{u}) - \vm{\varepsilon}_0) + \vm{\sigma}_0   &&  {\rm in }\; \Omega           
  \label{Eq:ConstitutiveRel} 
\end{align}
where  $\vm{L}$ is the differential operator for linear elasticity, and $\vm{G}$ is the projection operator that projects the stress field into tractions over any boundary, with $\vm{n}$ the unit normal to $\Gamma_N$, such that 
\begin{equation}
 \vm{L}^T = 
\begin{bmatrix}
  \partial/\partial x & 0                   & \partial/\partial y \\ 
  0                   & \partial/\partial y & \partial/\partial x
\end{bmatrix},
\qquad
 \vm{G} = 
\begin{bmatrix}
  n_x & 0 & n_y \\ 0 & n_y & n_x
\end{bmatrix},
\end{equation}
$\vm{D}$ is the matrix of the linear constitutive relation for stress and strain. We consider an homogeneous  Dirichlet boundary condition in (\ref{Eq:Dirichlet}) for simplicity. 

The problem expressed in its variational form is written as: 
\begin{equation}\label{Eq:WeakForm} 
 \begin{aligned}
  &\text{Find } \vm{u} \in V \text{ such that }  \forall \vm{v} \in V = \{\vm{v} \;|\; \vm{v} \in  [H^1(\Omega)]^2 , \vm{v}|_{\Gamma_D} = \vm{0} \}: \\ &
  \int _{\Omega}  
  {\vm{\varepsilon} (\vm{u})^{T} \vm{D} \vm{\varepsilon}(\vm{v}) \dd \Omega} =
\int _{\Omega} {\vm{v}^{T} \vm{b}   \dd \Omega} + 
\int _{\Gamma_N} {\vm{v}^{T} \vm{t}  \dd\Gamma} + 
\int _{\Omega} {\vm{\varepsilon} (\vm{v})^{T}\vm{D} \vm{\varepsilon}_0   \dd \Omega} -
\int _{\Omega} {\vm{\varepsilon}^{T}(\vm{v}) \vm{\sigma}_0    \dd \Omega}
 \end{aligned}
\end{equation}

% The symmetric and bilinear form $a: V \times V \rightarrow \mathbb{R}$ and the continuous linear form $l:V \rightarrow \mathbb{R}$ are defined in vectorial form by
% %
% \begin{align} \label{Eq:WeakForm1}
% a(\vm{u},\vm{v})  & := \int _{\Omega}  \vm{\sigma}^{T} (\vm{u})  \vm{\varepsilon}(\vm{v}) d \Omega = 
% \int _{\Omega}  \vm{\sigma}(\vm{u})^{T}  \vm{D}^{-1} \vm{\sigma} (\vm{v}) d \Omega \\
% l(\vm{v})&:= 
% \int _{\Omega} \vm{v}^{T} \vm{b}   \dd \Omega + \int _{\Gamma_N} \vm{v}^{T} \vm{t}  \dd\Gamma +
% \int _{\Omega} \vm{\sigma}^{T} (\vm{v}) \vm{\varepsilon}_0   \dd \Omega -
% \int _{\Omega} \vm{\varepsilon}^{T}(\vm{v}) \vm{\sigma}_0    \dd \Omega,\label{Eq:WeakForm2} 
% \end{align}
%  
% \noindent where $ \vm{\sigma}$  represents the stresses,  $\vm{\varepsilon}$ are the strains,  $\vm{D} $ is the elasticity matrix of the constitutive relation $\vm{\sigma}= \vm{D} \vm{\varepsilon} $ and $\vm{\sigma}_0$, $\vm{\varepsilon}_0$ the initial stresses and strains.

Let us consider a finite element approximation of $\vm{u}$ denoted as $\vm{u}^{h}$. In the XFEM formulation \cite{moesdolbow1999}, the approximation is usually enriched with two types of enrichment functions by means of the partition of unity: (i) a Heaviside function $H$ to describe the discontinuity of the displacement field along the crack, in the set of nodes $I^{crack}$ whose support is intersected by the crack and (ii) a set of branch functions $F_{\ell }$ to represent the asymptotic behaviour of the stress field near the crack tip, in the set of nodes $I^{tip}$ whose support contains the singularity. The XFEM displacement interpolation in a 2D model reads:
\begin{equation} \label{Eq:uXFEM} 
\vm{u}^{h}(\vm{x}) =\sum _{i\in I} N_{i}(\vm{x}) \vm{a}_{i}  +\sum _{i\in I^{crack}}N_{i}(\vm{x}) H(\vm{x})\vm{b}_{i}  +\sum _{i\in I^{tip}}N_{i}(\vm{x}) \left(\sum _{\ell =1}^{4}F_{\ell } (\vm{x})\vm{c}_{i}^{\ell }  \right)  
\end{equation}

\noindent where $N_i$ denotes the classical shape functions associated with node $i$ and $\vm{a}$, $\vm{b}$, $\vm{c}$ are the unknown coefficients. The $F_{\ell }$ functions used in this paper for the 2D case are \cite{moesdolbow1999}:
\begin{equation} 
\left\{F_{\ell } \left(r,\phi \right)\right\}\equiv \sqrt{r} \left\{\sin \frac{\phi }{2} ,\cos \frac{\phi }{2} ,\sin \frac{\phi }{2} \sin \phi ,\cos \frac{\phi }{2} \sin \phi \right\}
\end{equation}

% \begin{figure}[!ht]
%     \centering
%     \includegraphics{2_NodalSets} 
%     \caption{Classification of nodes in XFEM. Fixed enrichment area of radius $r_e$}
%     \label{fig:NodalSets}
% \end{figure}

Considering the enriched finite-dimensional subspace $V^{h} \subset V$ spanned by locally supported finite element shape functions, we solve for a discrete solution $\vm{u}^h \in V^{h}$ of the variational problem in (\ref{Eq:WeakForm}) such that $\forall \vm{v} \in V^{h}$:
%
% \begin{equation} \label{Eq:DiscreteWeakForm}  
%  \qquad a(\vm{u}^{h},\vm{v}) = l(\vm{v})     
% \end{equation} 

\begin{multline} \label{Eq:DiscreteWeakForm}  
\int _{\Omega}  {\vm{\varepsilon} (\vm{u}^h)^{T} \vm{D} \vm{\varepsilon}(\vm{v}) \dd \Omega} = \int _{\Omega}  {\vm{\sigma}^{T}(\vm{u}^h)  \vm{D}^{-1} \vm{\sigma} (\vm{v}) \dd \Omega} = \\
\int _{\Omega}   {\vm{v}^{T} \vm{b}   \dd \Omega} + 
\int _{\Gamma_N} { \vm{v}^{T} \vm{t}  \dd\Gamma} + 
\int _{\Omega}   {\vm{\varepsilon}(\vm{v})^{T} \vm{D} \vm{\varepsilon}_0   \dd \Omega} -
\int _{\Omega}   {\vm{\varepsilon}(\vm{v})^{T} \vm{\sigma}_0    \dd \Omega }
\end{multline}

%==============================================================================
\section{Error estimates in energy norm}
\label{sec:error}
%==============================================================================

\subsection{Zienkiewicz Zhu error estimate}

The discretisation error is defined as $\vm{e} := \vm{u}-\vm{u}^h$, in the absence of other types of errors. To quantify the error introduced by the discretisation a common approach is to use the energy norm of $\vm{e}$ defined as:
\begin{equation} \label{Eq:errorNorm} 
\norm{\vm{e}}^{2} = \int _{\Omega} \vm{\varepsilon}(\vm{e})^{T} \vm{D}  \vm{\varepsilon}(\vm{e})  \dd\Omega .  
\end{equation} 
Using the constitutive relation and introducing the error in the stress field $\vm{e}_\sigma := \vm{\sigma}- \vm{\sigma}^h $, where   $\vm{\sigma}^h = \vm{D}\left( \vm{\varepsilon}(\vm{u}^h)-\vm{\varepsilon}_0\right) + \vm{\sigma}_0$ is the finite element stress field, the previous expression can be written as
\begin{equation}
\norm{\vm{e}}^{2} =  \int _{\Omega}\vm{e}_\sigma ^{T} \vm{D}^{-1} \vm{e}_\sigma  \dd\Omega
 \end{equation}
Whereas the exact field $\vm{u}$ is in general unknown, it is possible to obtain an estimate of the error by means of the approximation introduced in \cite{zienkiewiczzhu1987} in the context of FE elasticity problems 
\begin{equation} \label{Eq:ZZ-estimator} 
\norm{\vm{e}}^{2} \approx  \int _{\Omega}\left( \vm{e}_\sigma^* \right)^{T} \vm{D}^{-1} \left( \vm{e}_\sigma^* \right) \dd\Omega ,  
\end{equation}
where $\vm{e}_\sigma^*$ is the approximated stress error defined by $\vm{e}_\sigma^* := \vm{\sigma}^* - \vm{\sigma}^h$, being $\vm{\sigma}^*$ the recovered stress field. Local element contributions are also obtained from (\ref{Eq:ZZ-estimator}) considering the domain of the element $\Omega_e$.

%==============================================================================
\subsection{Recovery technique}
\label{sec:recovery}
%==============================================================================

The accuracy of the Zienkiewicz-Zhu error estimator shown in (\ref{Eq:ZZ-estimator}) depends on the quality of the recovered field $\vm{\sigma}^*$. In this work we consider the SPR-CX recovery technique, which is an enhancement of the error estimator introduced in \cite{diezrodenas2007}, to recover the solutions for the primal and dual problems. The technique incorporates the ideas in \cite{rodenastur2007} to guarantee locally on patches the exact satisfaction of the equilibrium equations, and the extension in \cite{rodenasgonzalez2008} to singular problems. 

Let us define the field $\vm{\sigma}^-$ such that we subtract the initial stress and strain from the field $\vm{\sigma}$:
\begin{equation}\label{Eq:subtractIni}
 \vm{\sigma}^- = \vm{\sigma}-\vm{\sigma}_0 + \vm{D}\vm{\varepsilon_0},
\end{equation}
and perform the recovery on $\vm{\sigma}^-$. Then, the recovered field is
\begin{equation}
\vm{\sigma}^* = (\vm{\sigma}^{-})^* +\vm{\sigma}_0 - \vm{D}\vm{\varepsilon_0},
\end{equation}
where $(\vm{\sigma}^-)^*$ is the smoothed field that corresponds to $\vm{\sigma}^-$.

In the SPR-CX technique, as in the original SPR technique, we  define a patch $\mathcal{P}^{(J)}$ as the set of elements connected to a vertex node $J$. On each patch, a polynomial expansion for each one of the components of the recovered stress field is expressed in the form:
\begin{equation}
\hat{\sigma}_{k} ^{*} (\vm{x})  = \vm{p}(\vm{x}) \vm{a}_k \quad k=xx,yy,xy
\end{equation}
where $\vm{p}$ represents a polynomial basis and $\vm{a}_k$ are unknown coefficients. Usually, the polynomial basis is chosen equal to the finite element basis for the displacements. A least squares approximation to the values of FE stresses evaluated at the integration points of the elements within the patch, $\vm{x}_G \in \mathcal{P}^{(J)}$, is used to evaluate the coefficients $\vm{a}_k$. 

% \begin{align}
% \vm{p} (\vm{x}) &= \{1 \; x \; y \; x^2 \; xy \; y^2 \; \ldots \} \\
% \vm{a}_i (\vm{x}) &= \{ a_{0_i}(\vm{x}) \; a_{1_i}(\vm{x}) \;a_{2_i}(\vm{x}) \;a_{3_i}(\vm{x}) \;a_{4_i}(\vm{x}) \; a_{5_i}(\vm{x}) \; \ldots \}
% \end{align}

For the 2D case, the linear system of equations to evaluate the recovered stress field coupling the three stress components reads:
\begin{equation}
 \hat{\vm{\sigma}}^{*} (\vm{x}) = 
\begin{Bmatrix}
 \hat{\sigma}_{xx}^{*}(\vm{x})\\
 \hat{\sigma}_{yy}^{*}(\vm{x})\\  
 \hat{\sigma}_{xy}^{*}(\vm{x})                                               
\end{Bmatrix} = 
\vm{P}(\vm{x}) \vm{A} =
\begin{bmatrix}
 \vm{p}(\vm{x}) & \vm{0} & \vm{0} \\
 \vm{0} & \vm{p}(\vm{x}) & \vm{0} \\
 \vm{0} & \vm{0} & \vm{p}(\vm{x}) 
\end{bmatrix}
\begin{Bmatrix}
 \vm{a}_{xx}\\
 \vm{a}_{yy}\\  
 \vm{a}_{xy}  
\end{Bmatrix}
\end{equation}

In the basic SPR, we obtain the coefficients $\vm{A}$ from the minimisation of the functional 
\begin{equation}\label{Eq:MinFunctional}
 \mathcal{F}^{(J)}(A) = 
  \int_{\mathcal{P}^{(J)}} (\vm{P}  \vm{A} - \vm{\sigma}^{-h} )^2 \dd\Omega
\end{equation}
where $\vm{\sigma}^{-h} = \vm{D}\vm{\varepsilon}(\vm{u}^h)$.

The continuity of the recovered field is obtained by using a partition of unity procedure \cite{blackerbelytschko1994} to weight the stress fields obtained from the patches formed at the vertex nodes of the element. The field $\vm{\sigma}^*$ is interpolated using linear shape functions $N^{(J)}$ associated with the $n_v$ vertex nodes such that
\begin{equation}\label{Eq:conjoint_polynomials} 
 \vm{\sigma}^*(\vm{x}) = \sum_{J=1}^{n_v} N^{(J)}(\vm{x})  \hat{\vm{\sigma}}^{*(J)} (\vm{x}) - \vm{D}\vm{\varepsilon}_0(\vm{x}) + \vm{\sigma}_0(\vm{x}).
\end{equation}
Note that in \eqref{Eq:conjoint_polynomials} we add back the contribution of the initial stresses and strains subtracted in \eqref{Eq:subtractIni}.

\subsection{Equilibrium conditions}

Constraint equations are introduced via Lagrange multipliers into the functional defined in \eqref{Eq:MinFunctional} on each patch, in order to enforce the satisfaction of the:
\begin{itemize}
    \item Internal equilibrium equation: The constraint equation for the internal equilibrium in the patch is defined as:
    \begin{equation}
     \forall \vm{x_j} \in \mathcal{P}^{(J)} \qquad   \vm{L}^T  \hat{\vm{\sigma}}^{*(J)}(\vm{x}_j)  + \vm{L}^T
     (\vm{\sigma}_0(\vm{x}_j) - \vm{D}\vm{\varepsilon}_0(\vm{x}_j)) +  \hat{\vm{b}}(\vm{x}_j)  := \vm{c}^{\rm int}(\vm{x}_j) = 0    
    \end{equation}
    where $\hat{\vm{b}}(\vm{x})$ is a polynomial least squares fit of degree $p-1$ to the actual body forces $\vm{b}(\vm{x})$, being $p$ the degree of the recovered stress field $\hat{\vm{\sigma}}^{*(J)}$. We enforce $\vm{c}^{\rm int}(\vm{x}_j)$ at a sufficient number of $j$ non-aligned points ($nie$) to guarantee the exact representation of $\hat{\vm{b}}(\vm{x})$.
    
    \item Boundary equilibrium equations: We use a point collocation approach to impose the satisfaction of a polynomial approximation to the tractions along the Neumann boundary intersecting the patch. The constraint equation reads 
    \begin{equation}
    \forall \vm{x}_j \in \Gamma_N \cap \mathcal{P}^{(J)}\qquad \vm{G}  \hat{\vm{\sigma}}^{*(J)}(\vm{x}_j) + \vm{G} \vm{L}^T
     (\vm{\sigma}_0(\vm{x}_j) - \vm{D}\vm{\varepsilon}_0(\vm{x}_j)) - \vm{t}(\vm{x}_j) := \vm{c}^{\rm ext}(\vm{x}_j)  = 0      
    \end{equation}
    We enforce $\vm{c}^{\rm ext}(\vm{x}_j)$ in $nbe = p+1$ points along the part of the boundary crossing the patch. In the case that more than one boundary intersects the patch, only one curve is considered in order to avoid over-constraining.  
    
    \item Compatibility equations: $\vm{c}^{\rm cmp}(\vm{x}_j)$  is only imposed in the case that $p \geq 2$ in a sufficient number of non-aligned points. $\hat{\vm{\sigma}}^*$ directly satisfies $\vm{c}^{cmp}$ for $p=1$. 
\end{itemize}

Thus, the Lagrange functional enforcing the constraint equations for a patch $\mathcal{P}^{(J)}$ can be written as
\begin{equation}\label{Eq:LagrangeFunctional}
 \mathcal{L}^{(J)} (\vm{A},\vm{\lambda}) = 
 \mathcal{F}^{(J)}(\vm{A}) +
 \sum_{i=1}^{nie}\lambda_i^{\rm int}\left( \vm{c}^{\rm int}(\vm{x}_i)\right) + 
 \sum_{j=1}^{nbe}\lambda_j^{\rm ext}\left(\vm{c}^{\rm ext}(\vm{x}_j)\right) +
\sum_{k=1}^{nc}\lambda_k^{\rm cmp}\left(\vm{c}^{\rm cmp}(\vm{x}_k)\right).
\end{equation}

\begin{figure}[htb!]
    \centering
    \includegraphics{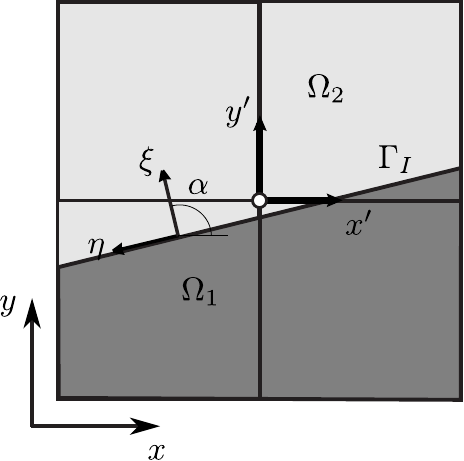}
    \caption{Equilibrium conditions along internal boundaries.}
    \label{fig:InternalBoundary}
\end{figure}

Optimizing functional \eqref{Eq:LagrangeFunctional} we obtain a linear system of equations to evaluate the coefficients $\vm{A}$. To enforce equilibrium conditions along internal boundaries (e.g. bimaterial problems, problems with zones subjected to different body forces, etc.), we consider different polynomial expansions on each side of the boundary and enforce the statical admissibility condition imposing equilibrium along this boundary. Suppose that we have a patch intersected by $\Gamma_I$ such that $\Omega_e = \Omega_{1,e} \cup \Omega_{2,e}$ for intersected elements, as shown in Figure \ref{fig:InternalBoundary}. To enforce equilibrium conditions along $\Gamma_I$ we define the stresses $\hat{\vm{\sigma}}^*_{\Omega_1} $, $\hat{\vm{\sigma}}^*_{\Omega_2} $ at each side of the internal boundary. Then, the boundary equilibrium along $\Gamma_I$ given the prescribed tractions $ \vm{t}_{\Gamma_I}= [t_x \; t_y]^T$ is:
\begin{equation}
\vm{G}(\hat{\vm{\sigma}}^*_{\Omega_1}|_{\Gamma_I} -\hat{\vm{\sigma}}^*_{\Omega_2}|_{\Gamma_I})  = \vm{t}_{\Gamma_I}.
\label{eq:ContIntEqDispl}
\end{equation}
%
% and $\vm{G}$ is the matrix form of Cauchy's law considering the unit normal $\vm{n}$ to $\Gamma_I$ such that
% \begin{equation}
%  \vm{G} = 
% \begin{bmatrix}
%   n_x & 0 & n_y \\ 0 & n_y & n_x
% \end{bmatrix}
% \end{equation}
% 

The same procedure can be used for patches intersected by the crack. In this case, we could consider the traction-free condition along the crack faces or define a different prescribed condition depending on the configuration.

After evaluating the equilibrated recovered fields on each patch $\hat{\vm{\sigma}}^{*(J)}$, we use \eqref{Eq:conjoint_polynomials} to obtain a continuous field. This process introduces a lack of equilibrium $\vm{s} = \sum _{J=1} ^{n_v} \nabla N^{(J)} \vm{\sigma}^{*(J)}$ when evaluating the divergence of the internal equilibrium equation, as explained in \cite{diezrodenas2007, rodenasgonzalez2010}. 

\subsection{Singular fields}

Different techniques have been used to account for the singular part during the recovery process \cite{rodenasgonzalez2008, bordasduflot2007}. Here, following the ideas in \cite{rodenasgonzalez2008}, for singular problems the exact stress field $\vm{\sigma}$  is decomposed into two stress fields, a smooth field $\vm{\sigma}_{\rm smo}$ and a singular field $\vm{\sigma}_{\rm sing}$:
\begin{equation} \label{Eq:splitting} 
\vm{\sigma} =   \vm{\sigma}_{\rm smo} + \vm{\sigma}_{\rm sing}.
\end{equation}

Then, the recovered field $\hat{\vm{\sigma}}^*$ required to compute the error estimate given in (\ref{Eq:ZZ-estimator}) can be expressed as the contribution of two recovered stress fields, one smooth $\hat{\vm{\sigma}}^*_{\rm smo}$ and one singular $\hat{\vm{\sigma}}^*_{\rm sing}$:
\begin{equation} \label{Eq:recovered_splitting} 
\hat{\vm{\sigma}}^* =   \hat{\vm{\sigma}}^*_{\rm smo} + \hat{\vm{\sigma}}^*_{\rm sing}.
\end{equation}

For the recovery of the singular part, the expressions which describe the asymptotic fields near the crack tip are used. To evaluate $\hat{\vm{\sigma}}^*_{\rm sing}$ we first obtain estimated values of the stress intensity factors $K_{\rm I}$ and $K_{\rm II}$ using a domain integral method based on extraction functions \cite{szabobabuska1991, ginertur2009b}. Notice that the recovered part $\hat{\vm{\sigma}}^*_{\rm sing}$ is an equilibrated field as it satisfies the equilibrium equations.

Once the field $\hat{\vm{\sigma}}^*_{\rm sing}$ has been evaluated, an FE-type approximation (discontinuous) to the smooth part $\hat{\vm{\sigma}}^h_{\rm smo}$ can be obtained subtracting $\hat{\vm{\sigma}}^*_{\rm sing}$ from the raw FE field:
\begin{equation} \label{Eq:tens_h_smooth} 
\hat{\vm{\sigma}}^h_{\rm smo} =   \hat{\vm{\sigma}}^h - \hat{\vm{\sigma}}^*_{\rm sing}.
\end{equation}

Then, the field $\hat{\vm{\sigma}}^*_{\rm smo}$ is evaluated applying the enhancements of the SPR technique previously described, i.e. satisfaction of equilibrium and compatibility equations at each patch. Note that as both $\hat{\vm{\sigma}}^*_{\rm smo}$  and $\hat{\vm{\sigma}}^*_{\rm sing}$ satisfy the equilibrium equations, $\hat{\vm{\sigma}}^*$ also satisfies equilibrium at each patch.

\section{Error in quantities of interest}
\label{sec:errorQoI}

\subsection{Exact error representation and auxiliary problem}
\label{sec:DualApproach}

The goal of many numerical computations is to control a specific design parameter, thus, it results natural to formulate the error in terms of such quantity. For this purpose, error estimators measured in the energy norm might be utilised to estimate the error in a particular quantity of interest \cite{ainsworthoden2000}. In this section we show how the ZZ estimate with the SPR-CX recovery may be used to evaluate the error in quantities of interest.

A common approach to evaluate the error in QoI involves the use of duality techniques which solve two different problems. A \emph{primal problem}, which is the problem at hand as shown in (\ref{Eq:WeakForm}), and a \emph{dual problem} used to extract information on the QoI. Thus, we shall explain the 
formulation of the dual problem.

Consider the \emph{primal} problem given in (\ref{Eq:WeakForm}) and its approximate finite element solution  $\vm{u}^h \in V^h \subset V$. Let $Q: V \rightarrow \mathbb{R}$ be a bounded linear functional representing some quantity of interest, acting on the space $V$ of admissible functions for the problem at hand. We are interested in estimating the error in the functional $Q (\vm{u})$ when calculated using the value of the approximate solution $\vm{u}^h$:
\begin{equation}
Q (\vm{u}) - Q (\vm{u}^h) = Q (\vm{u} - \vm{u}^h) = Q (\vm{e}) 
\end{equation}
%
%$Q(\vm{v})$ may be interpreted as the work associated with a displacement field $\vm{v}$ and a distribution of forces specific to each type of QoI. For $\vm{v} = \vm{u}$ in $Q (\vm{v})$, this force distribution can be used to extract information concerning the QoI associated with the solution of the problem in (\ref{Eq:WeakForm}).

To evaluate $Q (\vm{e})$  the standard procedure is to solve the auxiliary or \emph{dual} problem
%
% \begin{equation} 
%   a(\vm{v},\vm{w}_Q)=Q(\vm{v}). 
% \end{equation}

\begin{equation}\label{Eq:dual}
\begin{array}{l}
\displaystyle \text{Find } \tilde{\vm{u}} \in V  \text{ such that }  \forall \vm{v} \in V,  \\
\displaystyle \int _{\Omega}  \vm{\varepsilon}(\vm{v})^{T}  \vm{D}^{} \vm{\varepsilon} (\tilde{\vm{u}}) \dd \Omega = Q(\vm{v}),
\end{array}
\end{equation}
which can be seen as the variational form of an auxiliary mechanical problem used to extract information of the QoI. The dual displacement field $\tilde{\vm{u}} \in V$ vanishes over $\Gamma_D$. Test function $\vm{v}$ is a virtual displacement. Field $\tilde{\vm{\sigma}} = \vm{D}(\vm{\varepsilon}(\tilde{\vm{u}}) - \tilde{\vm{\varepsilon}}_0 ) +\tilde{\vm{\sigma}}_0$, where $\tilde{\vm{\sigma}}_0$ and $ \tilde{\vm{\varepsilon}}_0$ are known initial stress and strain, can be interpreted as a mechanical stress field. The left-hand side of \eqref{Eq:dual} is the work of internal forces of the auxiliary mechanical problem and $Q(\vm{v})$ is the work of an abstract external load.

We consider the same finite element space used in the primal problem to look for an approximation of $\tilde{\vm{u}} \in V$ such that the problem is
\begin{equation}\label{Eq:dualapprox}
\begin{array}{l}
\displaystyle \text{Find } \tilde{\vm{u}}^h \in V^h  \text{ such that }  \forall \vm{v} \in V^h,  \\
\displaystyle \int _{\Omega}  \vm{\varepsilon}(\vm{v})^{T}  \vm{D}^{} \vm{\varepsilon} (\tilde{\vm{u}}^h) \dd \Omega = Q(\vm{v}).
\end{array}
\end{equation}
To obtain an exact representation for the error $Q(\vm{e})$ in terms of the solution of the dual problem we substitute $\vm{v}=\vm{e}$ in (\ref{Eq:dual}) and, considering the Galerkin orthogonality, for all $\tilde{\vm{u}}^h \in V^h$:
\begin{align}
 \label{Eq:GalerkinOrtho}
Q(\vm{e}) = \int _{\Omega}  \vm{\varepsilon}(\vm{e})^{T}  \vm{D}^{} \vm{\varepsilon} (\tilde{\vm{e}}) \dd \Omega
\end{align}
where $\tilde{\vm{e}} := \tilde{\vm{u}} - \tilde{\vm{u}}^h$ is the discretisation error of the dual problem \eqref{Eq:dual}. We can obtain an expression in terms of the mechanical stresses using the constitutive relation:
\begin{equation}\label{Eq:errorQoI}
Q(\vm{e}) =  \int_\Omega \vm{e}_{\vm{\sigma}}^{T}  \vm{D}^{-1} \tilde{\vm{e}}_{\vm{\sigma}}  \dd \Omega
\end{equation}
where $\tilde{\vm{e}}_\sigma := \tilde{\vm{\sigma}} - \tilde{\vm{\sigma}}^{h}$ is the stress error of the dual problem and $\tilde{\vm{\sigma}}^{h} = \vm{D}  (\vm{\varepsilon} (\tilde{\vm{u}}^h) - \tilde{\vm{\varepsilon}}_0) + \tilde{\vm{\sigma}}_{0}$ the finite element stress field. 
% An exact representation for the error $Q(\vm{e})$ in terms of the solution of the dual problem can be simply obtained by substituting $\vm{v}=\vm{e}$ in (\ref{Eq:dual}) and remarking that for all $\vm{w}_Q^h \in V^h$, due to the Galerkin orthogonality, $a(\vm{e},\vm{w}_Q^h)=0$  such that
% %
% \begin{equation}
% Q (\vm{e}) = a (\vm{e}, \vm{w}_Q) = a (\vm{e}, \vm{w}_Q)  -  a (\vm{e}, \vm{w}_Q^h) = a (\vm{e}, \vm{w}_Q  -  \vm{w}_Q^h) = a (\vm{e}, \vm{e}_Q)
% \end{equation}

\subsection{Smoothing-based error estimate}

The error in the QoI in (\ref{Eq:errorQoI}) is related to the errors in the FE approximations $\vm{u}^h$ and $\tilde{\vm{u}}^h$. Thus, we can select from the set of available procedures to estimate the error in the energy norm a technique to obtain estimates of the error in the QoI. Considering expressions (\ref{Eq:ZZ-estimator}) and (\ref{Eq:errorQoI}) we can derive an estimate for the error in the QoI which reads
\begin{equation}\label{Eq:errestQoI}
Q(\vm{e}) \approx \mathcal{E} = \int_\Omega {(\vm{e}_{\vm{\sigma}}^*)}^{T}  \vm{D}^{-1} {(\tilde{\vm{e}}_{\vm{\sigma}}^*)}  \dd \Omega
\end{equation}
where the approximate dual error is $\tilde{\vm{e}}_\sigma^* = \tilde{\vm{\sigma}}^{*} - \tilde{\vm{\sigma}}^{h}$ and $\tilde{\vm{\sigma}}^*$ is the recovered auxiliary stress field. Here, we expect to have a sharp estimate of the error in the QoI if the recovered stress fields are accurate approximations to their exact counterparts.

The recovered stress fields can be computed in many ways, for example, by using the SPR technique as explained in \cite{zienkiewiczzhu1992}. To obtain accurate representations of the exact stress fields for the primal and dual solutions, we propose the use of the locally equilibrated recovery technique described in Section~\ref{sec:recovery}. This technique, which is an enhancement of the SPR, enforces the fulfilment of the internal and boundary equilibrium equations locally on patches. For problems with singularities the stress field is also decomposed into two parts: smooth and singular, which are separately recovered. 

Two remarks have to be made. First, the analytical expressions that define the tractions and body forces for the dual problem are obtained from the interpretation of the functional $Q$ in terms of tractions, body loads, initial stresses and strains, as seen in Section \ref{sec:DualProblem}. Second, to enforce equilibrium conditions during the recovery process along the boundary of the domain of interest (DoI) used to define the QoI, we consider it as an internal interface. We use different polynomial expansions on each side of the boundary and enforce statical admissibility of the normal and tangential stresses as previously explained in Section \ref{sec:recovery}

%==============================================================================
\subsection{Analytical definitions for the dual problem}
\label{sec:DualProblem}
%==============================================================================

The SPR-CX recovery requires that the mechanical equilibrium must be made explicit in order to recover the dual stress field. Thus, the right-hand side of \eqref{Eq:dual} is interpreted as the work of mechanical external forces, and the analytical expression of these forces is derived, depending on the quantity of interest: 
\begin{equation} \label{Eq:WeakFormDual}
\begin{aligned} 
&\text{Find } \tilde{\vm{u}} \in  V  \text{ such that }\forall \vm{v} \in V: \\
&\int _{\Omega}  {\vm{\varepsilon} (\vm{v})^{T} \vm{D}   \vm{\varepsilon}(\tilde{\vm{u}}) \dd \Omega} = Q(\vm{v})\\ 
&=\int _{\Omega} {\vm{v}^{T} \tilde{\vm{b}}   \dd \Omega} 
+\int _{\Gamma_N} {\vm{v}^{T} \tilde{\vm{t}}  \dd\Gamma} 
+\int _{\Omega}   { \vm{\varepsilon}(\vm{v})^T \vm{D} \tilde{\vm{\varepsilon}}_0   \dd \Omega} 
-\int _{\Omega} {\vm{\varepsilon}(\vm{v})^{T} \tilde{\vm{\sigma}}_0    \dd \Omega}
\end{aligned} 
\end{equation}

The problem in \eqref{Eq:WeakFormDual} is solved using a FE approximation with test and trial functions in $V^h$. The finite element solution is denoted by $\tilde{\vm{u}}^h \in V^h$.

Such derivations were presented in \cite{rodenas2005,gonzalezrodenas2011b,verdugodiez2011}. Here, we only recall some of the results presented in these papers. Additionally, we provide the analytical expression of the dual load when the quantity of interest is the generalised stress intensity factor (GSIF).

\subsubsection{Mean strain in \texorpdfstring{$\Omega_I$}{Omegai}}
%-------------------------------------------

In this case we are interested in some combination of the components of the strain over a subdomain $\Omega_I$ such that the QoI is given by:

\begin{equation}\label{Eq:epsmean}
 Q(\vm{u}) 
 %= \bar{\varepsilon}  |_{\Omega_I} 
  = \frac{1}{|\Omega_I|} 
  \int_{\Omega_I} {\vm{c}_{\varepsilon}^T \vm{\varepsilon}(\vm{u}) \dd \Omega}  
  = \int_{\Omega_I} {\frac{\vm{c}_{\varepsilon}^T}{|\Omega_I|} \vm{\varepsilon}(\vm{u}) \dd \Omega}
\end{equation}
where $\vm{c}_{\varepsilon}$ is the extraction operator used to define the combination of strains under consideration. Thus, the term $\tilde{\vm{\sigma}}_0 = \vm{c}_{\varepsilon}^T/|\Omega_I|$ represents the vector of initial stresses that are used to define the auxiliary problem for this particular QoI.

%-------------------------------------------
\subsubsection{Mean stress value in \texorpdfstring{$\Omega_I$}{Omegai}}
%-------------------------------------------
\label{sec:MeanStressOmega}

Let us consider now as QoI the mean stress value 
%$\bar{\sigma} |_{\Omega_I}$ 
given by a combination $\vm{c}_{\sigma}$ of the stress components $
\vm{\sigma} = \vm{D} ( \vm{\varepsilon}(\vm{u}) - \vm{\varepsilon}_0) + \vm{\sigma}_0$ in a domain of interest which reads:
\begin{equation}
  Q(\textbf{u}) 
  %=   \bar{\sigma}  |_{\Omega_I} 
%   = \frac{1}{\left|\Omega_I\right|}\int_{\Omega_I} \sigma (\vm{u})  \dd \Omega
  = \frac{1}{|\Omega_I|}\int_{\Omega_I} \vm{c}^T_{\sigma} (\vm{D} ( \vm{\varepsilon}(\vm{u}) - \vm{\varepsilon}_0) + \vm{\sigma}_0)\dd \Omega .
  \label{eq:MoI_mStress}
\end{equation}

$Q$ is an affine functional. Let us define
\begin{equation}
 \tilde{Q}(\vm{v}) = \int _{\Omega}   \vm{c}^T_{\sigma} \vm{D} ( \vm{\varepsilon}(\vm{v}))\dd \Omega  
\end{equation}
for $\vm{v}$ an arbitrary vector of $H^1(\Omega)$. $\tilde{Q}$ is such that $\tilde{Q}(\vm{e}) = Q(\vm{e})$, so that by solving the dual problem
\begin{equation}
 \int _{\Omega}  \vm{\varepsilon} (\vm{v})^{T} \vm{D}   \vm{\varepsilon}(\tilde{\vm{u}}) \dd \Omega = 
\tilde{Q}(\vm{v})
\end{equation} 
for $\tilde{\vm{u}}$, the derivations of Section  \ref{sec:DualApproach} apply.

Similarly to the previous quantity, the right-hand side of the auxiliary problem is defined by the  term $\tilde{\vm{\varepsilon}}_0 = \vm{c}_{\sigma }^T/|\Omega_I|$, which represents in this case a vector of initial strains.

% %==============================================================================
% 
% \subsection{Mean stresses and strains in a domain of interest from the perspective of initial stresses and strains}
% 
% Previously, the expressions to calculate mean stresses or strains in a domain of interest $\Omega_i$ have been obtained using the divergence theorem. However, we can calculate these quantities of interest considering them as an initial stress load case in the dual problem. This can be proved immediately
% %
% \begin{equation}
%  Q(\vm{u}) = \bar{\sigma}_\alpha |_{\Omega_i} 
%   = \frac{1}{|\Omega_i|}\int_{\Omega_i} \vm{c}_{\sigma_\alpha}^T \vm{\sigma} \dd \Omega  
%   = \int_{\Omega_i} \frac{\vm{c}_{\sigma_\alpha}^T}{|\Omega_i|} \vm{\sigma} \dd \Omega
% \end{equation}
% 
% From \eqref{Eq:WeakForm} we can define $\vm{\varepsilon}_{0,d} =  {\vm{c}_{\sigma_\alpha}^T} /{|\Omega_i|}$  corresponding to the term of initial strains  that we need to apply in the dual problem to extract the value of $\bar{\sigma}_\alpha |_{\Omega_i}$. A similar formulation can be derived for the case of the mean strains in $\Omega_i$ such that $\vm{\sigma}_{0,d} =  {\vm{c}_{\varepsilon_\alpha}^T}/{|\Omega_i|}$.

%==============================================================================
\subsubsection{Generalised stress intensity factor}
\label{sec:SIF-QoI}
%==============================================================================

The generalised stress intensity factor (GSIF) $K$ is the characterizing parameter in singular problems as in the case of reentrant corners or in fracture mechanics. For that reason, it is important to evaluate error estimates considering this parameter as the quantity of interest. To evaluate the GSIF in XFEM approximations it is a common practice to use the interaction integral in its equivalent domain integral (EDI) form. There are different expressions already available to evaluate EDI integrals for singular problems. In this work, we are going to consider the method based on extraction functions, as shown in \cite{szabobabuska1991}, which is a generalisation of the interaction integral for the singular problem of a V-notch plate:
%szabobabuska1991 p219  
\begin{equation} \label{Eq:EDIIntegral} 
    Q(\vm{u}) = K = -\frac{1}{C} \int_{\Omega} \left( \sigma_{jk}^{} u_{k}^{\rm aux} -  \sigma_{jk}^{\rm aux}  u_{k}^{} \right) \frac{\partial{q}}{\partial{x_j}} \dd \Omega
\end{equation}
where $ u^{\rm aux}$, $\sigma^{\rm aux}$ are the auxiliary fields used to extract the GSIFs in mode I or mode II and $C$ is a constant that is dependent on the geometry and the loading mode. $q$ is an arbitrary $C^0$ function that defines the extraction zone $\Omega_I$ which takes the value of 1 at the singular point and 0 at the boundary $\Gamma$, $x_j$ refers to the local coordinate system defined at the singularity.

To formulate the dual problem, we assemble the vector of equivalent nodal forces corresponding to the volume loads in the domain of interest that represent the stress intensity factor. Consider the expanded expression for $K$ with three terms function of the primal stresses and two terms function of the primal displacements:
\begin{multline} 
   Q(\vm{u}) =  K =  
\int_{\Omega_I}  (\vm{\sigma}^{})^T \left(-\frac{1}{C}\right)
\begin{bmatrix}
u_{1}^{\rm aux} q_{,1} \\
u_{2}^{\rm aux} q_{,2}\\
u_{2}^{\rm aux}q_{,1} + u_{1}^{\rm aux}q_{,2}
\end{bmatrix}  - \\
 (\vm{u}_{}^{})^T \left(-\frac{1}{C}\right)
\begin{bmatrix}
\sigma_{11}^{\rm aux}q_{,1} + \sigma_{21}^{\rm aux} q_{,2} \\
\sigma_{12}^{\rm aux}q_{,1} + \sigma_{22}^{\rm aux} q_{,2}
\end{bmatrix} 
  \dd \Omega
\end{multline} 
which can be rewritten as a function of initial strains $\tilde{\vm{\varepsilon}}_0$ and body loads $\tilde{\vm{b}}$:
\begin{equation}
 Q(\vm{u}) =  K =  
\int_{\Omega_I}  \vm{\sigma}_{}^{}(\vm{u})^T \tilde{\vm{\varepsilon}}_0 + 
 (\vm{u}^{})^T \tilde{\vm{b}}    \dd \Omega.
\end{equation}

Thus, if we replace $\vm{u}$ with the vector of arbitrary displacements $\vm{v}$, the quantity of interest can be evaluated from
\begin{equation}
  Q(\vm{v}) =  
\int_{\Omega_I} \vm{\sigma}(\vm{v})^T \tilde{\vm{\varepsilon}}_{0} \dd \Omega + 
\int_{\Omega_I} \vm{v}_{}^{T} \tilde{\vm{b}} \dd \Omega .
\end{equation}

Hence, the initial strains and the body loads per unit volume that can be applied in the dual problem to extract the GSIF are defined as
\begin{equation} \label{Eq:IniStrainSIF}
 \tilde{\vm{\varepsilon}}_{0} = 
-\frac{1}{C}
\begin{bmatrix}
u_1^{\rm aux}q_{,1}\\
u_2^{\rm aux}q_{,2}\\ 
u_2^{\rm aux}q_{,1}+u_1^{\rm aux}q_{,2}
\end{bmatrix}  
\quad,\quad 
\tilde{\vm{b}} = \frac{1}{C} 
\begin{bmatrix}
\sigma_{11}^{\rm aux}q_{,1}+\sigma_{21}^{\rm aux}q_{,2}\\
\sigma_{12}^{\rm aux}q_{,1}+\sigma_{22}^{\rm aux}q_{,2}
\end{bmatrix}
\end{equation}

% \subsection{Local contributions}
% 
% The ZZ error estimate in quantities of interest given in equation \eqref{Eq:errestQoI} can be written in terms of local contributions, which proves useful for adaptivity purposes. 
% For a discretisation with $n_e$ elements, element $e$ occupying domain $\Omega_e$ such that $\Omega_e: \Omega= \bigcup_{n_e} \Omega_e$, we can write:
% \begin{equation}
%  \mathcal{E} = 
%  %\int_\Omega \vm{e}_\sigma^*  \vm{D}^{-1} \tilde{\vm{e}}_\sigma^*  \dd \Omega =
%  \sum_{n_e} \int _{\Omega_e}  
%  {\vm{e}_{\vm{\sigma}}^*} ^{T}  
%  \vm{D} ^{-1}
%  \tilde{\vm{e}}^*_{\vm{\sigma}}  
%  \dd \Omega
% \end{equation} 

%==============================================================================
\section{Numerical results}
\label{sec:results}
%==============================================================================
In this section we consider numerical examples for 2D problems with exact analytical solution to evaluate the performance of the technique presented above. For that purpose we define the effectivity index of the error estimator $\theta$ as: 
\begin{equation} \label{Eq:Effectivity}  
\theta  =\frac{\mathcal{E}}{Q(\vm{e})}   
\end{equation}  

\noindent where $Q(\vm{e})$ denotes the exact error in the quantity of interest, and  $\mathcal{E}$ represents the evaluated error estimate. We can also represent the effectivity in the QoI defined as
\begin{equation}
  \theta_{QoI} =\frac{Q(\vm{u}^h)+\mathcal{E}}{ Q(\vm{u})} 
\end{equation}

\noindent and the relative error in the QoI for the exact and estimated error 
\begin{equation}
 \eta^Q(\vm{e}) = \frac{\abs{Q( \vm{e} )}}{\abs{Q( \vm{u} )}} , \quad
 \eta^\mathcal{E} = \frac{\abs{\mathcal{E}}}{\abs{Q(\vm{u}^h)+\mathcal{E}}}
\end{equation}

\subsection{ Westergaard problem -- FEM solution. } 
 
\begin{figure}[ht]
    \centering
    \includegraphics{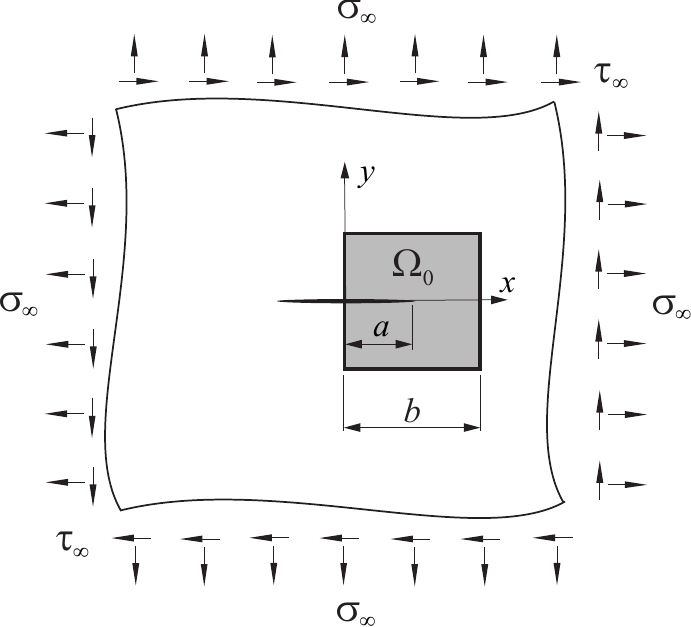}
    \caption{Westergaard problem. Infinite plate with a crack of length $2a$ under uniform tractions $\sigma_{\infty}$ (biaxial) and $\tau_{\infty}$. Finite portion of the domain $\Omega_0$, modelled with FE.}
    \label{fig:westergaard}
\end{figure}

Let us consider the Westergaard problem \cite{rodenasgonzalez2008, ginerfuenmayor2005} of linear elastic fracture mechanics for which the exact analytical solution is known. The Westergaard problem corresponds to an infinite plate loaded at infinity with biaxial tractions $\sigma_{x \infty}=\sigma_{y \infty}=\sigma_{\infty}$ and shear traction $\tau_{\infty}$, presenting a crack of length $2a$ as shown in Figure~\ref{fig:westergaard}. Combining the externally applied loads we can obtain different loading conditions: pure mode I, pure mode II or mixed mode.  
 
The numerical model corresponds to a finite portion of the domain ($a=5$ and $b=10$ in Figure~\ref{fig:westergaard}). The applied projected stresses for mode I are evaluated from the analytical Westergaard solution \cite{ginerfuenmayor2005}:  

\begin{equation} \label{Eq:WesterStressI} 
\begin{array}{r@{\hspace{1ex}}c@{\hspace{1ex}}l} 
{\sigma _{x}^{I} } (x,y) & {=} & {\displaystyle \frac{\sigma _{\infty } }{\sqrt{\left|t\right|} } \bigg[\left(x\cos \frac{\phi }{2} -y\sin \frac{\phi }{2} \right)+y\frac{a^{2} }{\left|t\right|^{2} } \left(m\sin \frac{\phi }{2} -n\cos \frac{\phi }{2} \right)\bigg]} \\ 
\noalign{\medskip}{\sigma _{y}^{I} }(x,y) & {=} & {\displaystyle \frac{\sigma _{\infty } }{\sqrt{\left|t\right|} } \bigg[\left(x\cos \frac{\phi }{2} -y\sin \frac{\phi }{2} \right)-y\frac{a^{2} }{\left|t\right|^{2} } \left(m\sin \frac{\phi }{2} -n\cos \frac{\phi }{2} \right)\bigg]} \\ 
\noalign{\medskip}{ \tau _{xy}^{I} }(x,y) & {=} & {\displaystyle y\frac{a^{2} \sigma _{\infty } }{\left|t\right|^{2} \sqrt{\left|t\right|} } \left(m\cos \frac{\phi }{2} +n\sin \frac{\phi }{2} \right)} \end{array}
\end{equation}
\noindent and for mode II:
\begin{equation} \label{Eq:WesterStressII} 
\begin{array}{r@{\hspace{1ex}}c@{\hspace{1ex}}l} 
{\sigma _{x}^{II}}(x,y)  & {=} & {\displaystyle \frac{\tau _{\infty } }{\sqrt{\left|t\right|} } \bigg[2\left(y\cos \frac{\phi }{2} +x\sin \frac{\phi }{2} \right)-y\frac{a^{2} }{\left|t\right|^{2} } \left(m\cos \frac{\phi }{2} +n\sin \frac{\phi }{2} \right)\bigg]} \\ 
\noalign{\medskip}{\sigma _{y}^{II}}(x,y)  & {=} & {\displaystyle y\frac{a^{2} \tau _{\infty } }{\left|t\right|^{2} \sqrt{\left|t\right|} } \left(m\cos \frac{\phi }{2} +n\sin \frac{\phi }{2} \right)} \\ \noalign{\medskip}{\tau _{xy}^{II}}(x,y)  & {=} & {\displaystyle \frac{\tau _{\infty } }{\sqrt{\left|t\right|} } \bigg[\left(x\cos \frac{\phi }{2} -y\sin \frac{\phi }{2} \right)+y\frac{a^{2} }{\left|t\right|^{2} } \left(m\sin \frac{\phi }{2} -n\cos \frac{\phi }{2} \right)\bigg]} \end{array}
\end{equation}

\noindent where the stress fields are expressed as a function of $x$ and $y$, with origin at the centre of the crack. The parameters $t$, $m$, $n$ and $\phi$ are defined as
\begin{equation}
\begin{split}t& =(x+iy)^{2} -a^{2} =(x^{2} -y^{2} -a^{2} )+i(2xy)=m+in \\  m & =\textrm{Re}(t) =\textrm{Re}(z^{2} -a^{2} )=x^{2} -y^{2} -a^{2} \\ n & =\textrm{Im}(t)=(z^{2} -a^{2} )=2xy \\ \phi & =\textrm{Arg} (\bar{t})=\textrm{Arg} (m-in) \qquad\textrm{with }\phi \in \left[-\pi ,\pi \right], \; i^2=-1 \end{split} 
\end{equation}

For the problem analysed, the exact value of the SIF is given by 
\begin{equation} \label{Eq:SIFWestergaard} 
K_{I,ex} =\sigma_{\infty } \sqrt{\pi a} \qquad \qquad K_{II,ex} =\tau_{\infty } \sqrt{\pi a}  
\end{equation}

Material parameters are Young's modulus $E = 10^7$ and Poisson's ratio $ \nu= 0.333$. We consider loading conditions in pure mode I with $\sigma_{\infty} =100$ and $\tau_{\infty}=0$, and pure  mode II with $\sigma_{\infty} =0$ and $\tau_{\infty}=100$. We assume plane strain conditions.

In the numerical analyses, we use a geometrical enrichment defined by a circular fixed enrichment area $B(x_0, r_e)$ with radius  $r_e = 2.5$, with its centre at the crack tip $x_0$ as proposed in \cite{bechetminnebo2005}. Bilinear elements are considered in the models, using a sequence of uniformly refined meshes. For the numerical integration of standard elements we use a $2\times2$ Gaussian quadrature rule. We use a $5\times5$ quasipolar integration in the subdomains of the element containing the crack tip \cite{bechetminnebo2005}. We do not consider correction for blending elements. Methods to address blending errors are proposed in \cite{chessawang2003, graciewang2008, fries2008, taranconvercher2009}. 

\begin{figure}[htb]
    \centering
    \includegraphics{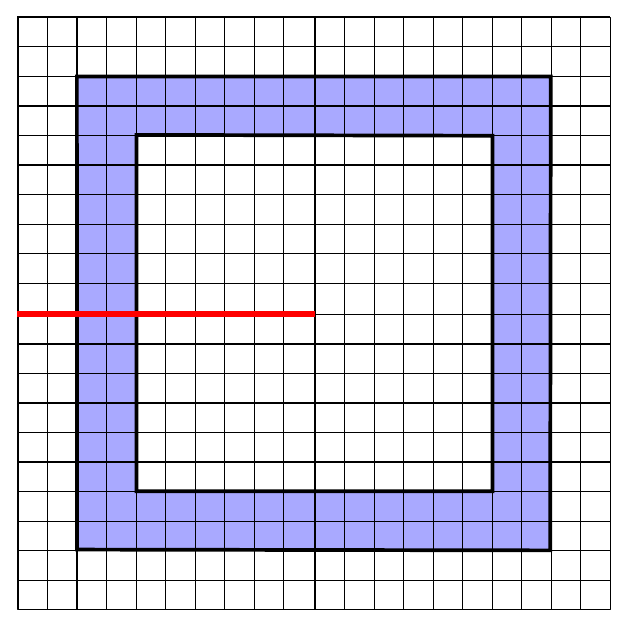}
    \caption{Domain of interest for the extraction of the stress intensity factor}
    \label{fig:DoISIF}
\end{figure}

To evaluate the stress intensity factor $K$ we use an EDI technique \cite{szabobabuska1991}. For the primal problem we consider a square plateau function $q$ centred at the crack tip, as shown in Figure \ref{fig:DoISIF}. $q=1$ for the domain defined by an inner square with side length 6 and $q=0$ for the part of the domain outside the outer square with side length 8, $q$ is interpolated in-between the two squares. This plateau function is also used to define the subdomain $\Omega_i$ when extracting the quantity of interest in the dual problem. As the dual problem is also a singular problem we have to evaluate a second stress intensity factor. In this case, we use a plateau function such that $q=1$ for all nodes inside a square with side length 4.9 and $q=0$ otherwise.
%TODO JJ: consider drawing the enrichment area + plateau for the dual problem. << Not very nice

\begin{figure}[htb]
    \centering
    \includegraphics{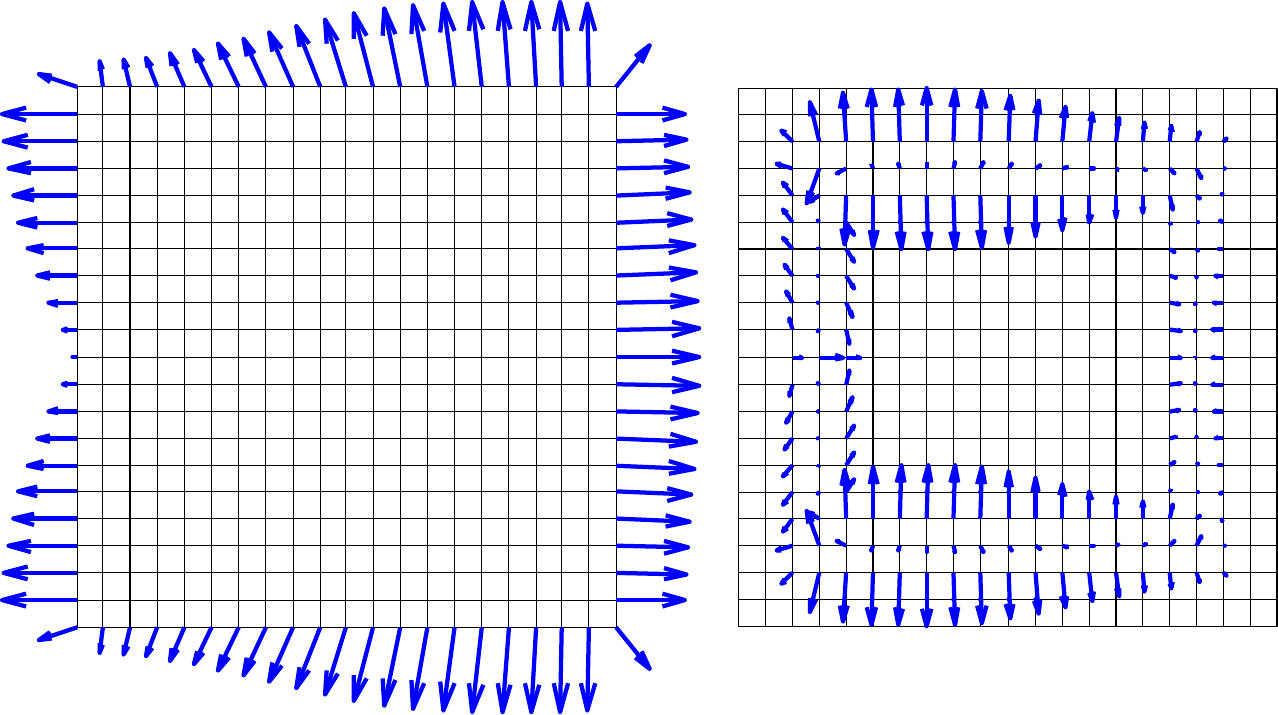}
    \includegraphics[scale=1.17]{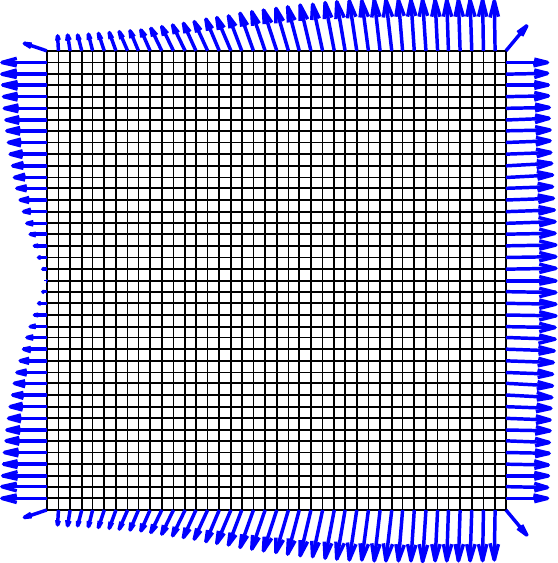}\hspace{20pt}
    \includegraphics[scale=1.17]{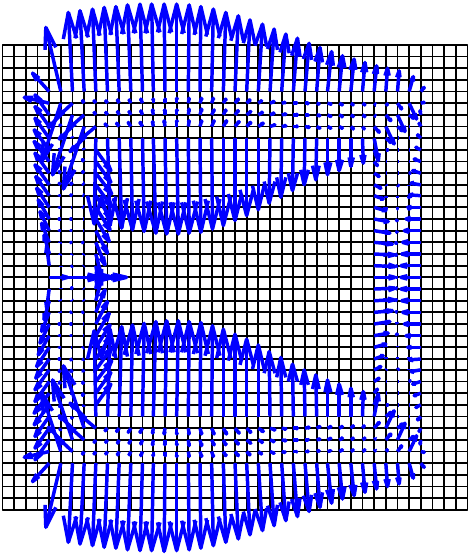}
    \caption{Equivalent forces at nodes for the primal (left) and dual (right) problems.}
    \label{fig:NodalForces}
\end{figure}

In Figure \ref{fig:NodalForces} we represent the equivalent nodal forces used to solve the primal and dual problems. For the dual problem the vector of forces is constructed using the discrete approximation of the dual function. The Dirichlet boundary constraints are the same for both models. For the dual problem, we can see that the forces are distributed in the nodes located in the domain of interest. For the recovery of the primal and dual fields we perform the splitting of stresses and enforce internal and boundary equilibrium and the compatibility equation. 
%To improve the computer efficiency we enforce internal equilibrium in the dual problem only outside the domain of interest, as this has been shown in \cite{gonzalezrodenas2013} to yield similar results with less computational effort. 

\begin{figure}[htb]
    \centering
    \includegraphics {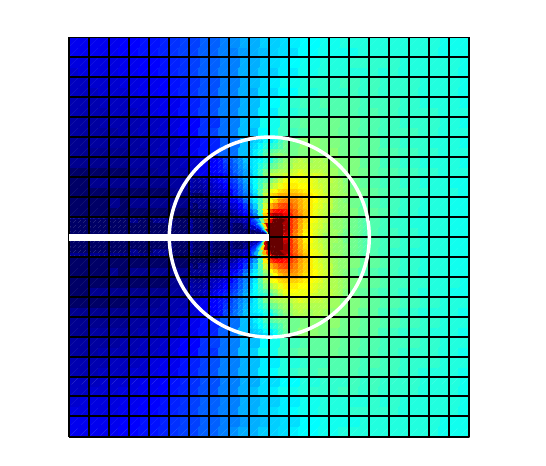}
    \includegraphics {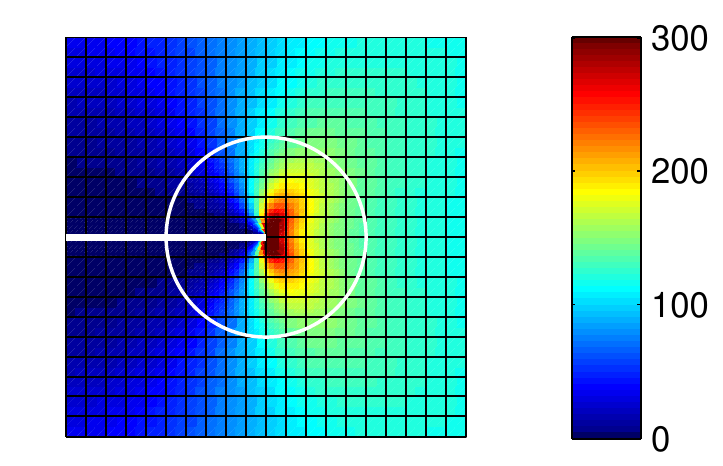}    
    \caption{FE (left) and recovered  (right) $\sigma_{yy}$ for the primal problem.}
    \label{fig:sigma_p}
\end{figure}

\begin{figure}[htb]
    \centering
    \includegraphics {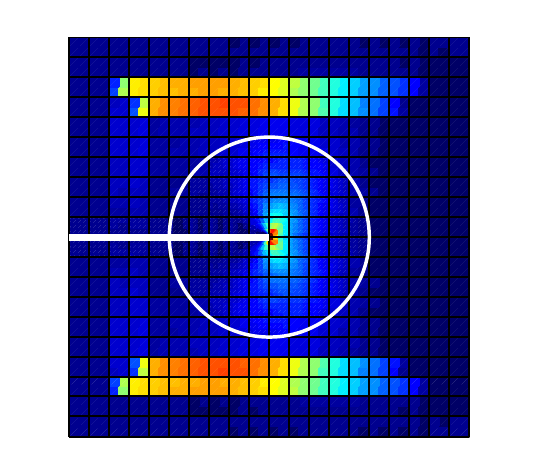}
    \includegraphics {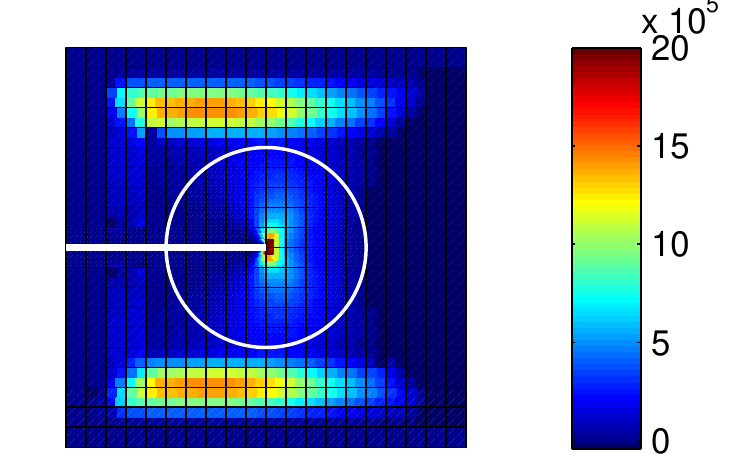}
    \caption{FE (left) and recovered fields (right) $\sigma_{yy}$ for the dual problem.}
    \label{fig:sigma_d}
\end{figure}

The \emph{yy}-component of the stress field for the raw FE and the recovered solutions is represented in Figure \ref{fig:sigma_p}. The enrichment area is indicated with a circle. In Figure \ref{fig:sigma_d} we show the same results for the dual problem. Notice how the recovery procedure smoothes the stresses along the interface of the domain of interest. As the dual problem is also characterised by the crack, we have to evaluate the corresponding stress intensity factor and perform the \emph{singular+smooth} decomposition of the stress field. 

%TODO the enrichment radius is 2, not 2.5 as in the figs.

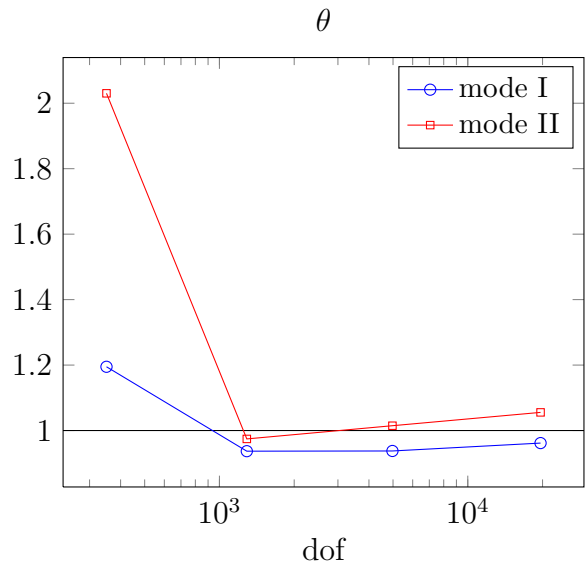
\begin{figure}[htb]
      \centering
 \begin{tikzpicture}%[scale=0.74]
    %\tikzset{/pgfplots/no markers}
    \begin{semilogxaxis}[
    %width=0.6\textwidth,
    title={$\theta$},
    xlabel={dof},
    legend style={cells={anchor=west}, font=\small},
    legend style={at={(0.98,0.98)},anchor=north east},
    cycle list name=ageplot
    ]
    \draw (axis cs: 0,1)--(axis cs: 1e8,1); 
    \addplot table[x=NDOF,y= Eff]{Data/CRACK_MI.dat};
    \addplot table[x=NDOF,y= Eff]{Data/CRACK_MII.dat};
    \legend{{mode I},{mode II}}
    \end{semilogxaxis}
\end{tikzpicture}
    \caption{Evolution of the effectivity index $\theta$ considering the SIF as quantity of interest under mode I and mode II loading conditions.}
\label{fig:theta}
\end{figure}

Figure \ref{fig:theta} shows the evolution of the effectivity  index $\theta$  as we increase the number of degrees of freedom (dof). We consider as quantities of interest the two GSIFs characterising two different loading conditions, i.e. mode I and mode II. We can see that for both quantities the error estimator yields effectivities close to the optimal value $\theta=1$. 

In Tables \ref{tab:modeI} and \ref{tab:modeII} we indicate the values for the estimated, $\mathcal{E}$, and exact, $Q(\vm{e})$, errors, the global effectivity index $\theta$ and the effectivity for the quantity of interest $\theta_{QoI}$ using the proposed recovery technique and the standard SPR (denoted with $\dagger$). The magnitude of the exact error is accurately captured by the estimated obtained with the SPR-CX, which is clearly reflected in the good effectivity index for both loading modes. As expected, the effectivity in the quantities of interest $\theta_{QoI}$ is highly accurate. For the SPR, although the value $\mathcal{E}$ decreases as we increase the number of dof, the error estimate is not as accurate as the estimate obtained with the SPR-CX and does not decrease as fast as the exact error, loosing asymptotic exactness. The SPR does not consider the splitting of the singular stresses, giving less accurate results close to the crack tip, and does not enforce equilibrium conditions in the primal and dual recovered fields, which results in a poorer description of the stresses close to the boundaries and the interface of the domain of interest. 

\begin{table}
\centering
\caption{Stress intensity factor $K_{\rm I}$ as QoI. $^{\dagger}$Results using the standard SPR recovery.}
\begin{filecontents*}{Data/CRACK_MI.csv}
dof,Qees,Qeex,theta,EffQoI,QeesSPR,QeSPR,Q(u^h)SPR,F_d*U_pSPR,thetaSPR,thetaQoISPR,etaSPR,etaesSPR,K_dSPR,E2SPR,theta_2SPR,E3SPR,theta_3SPR,E4SPR,theta_4SPR
351,2.5264002,2.1144455,1.1948288,1.0026054,13.7089424642199,2.11444545232055,155.999437556098,156.052512750098,6.48346943600494,1.07333003776324,1.33729272350358,8.07794079619323,1923083.73697604,14.3629114278795,6.79275571385234,20.7059309345672,9.79260586355774,36.5418920974749,17.2820216560191
1289,0.4821894,0.5146442,0.9369373,0.9997947,6.38796901145607,0.514644244577056,157.599238763842,157.614508816397,12.4123976489931,1.03714616740243,0.325489599512051,3.89540690284155,1710800.38299562,7.14614525112204,13.8856021930156,9.08747080635642,17.6577721447651,18.2300654192106,35.4226547975727
4973,0.1140116,0.1216039,0.9375659,0.999952,3.17557662748032,0.121603878078531,157.99227913034,157.996215472512,26.1141065372072,1.01931501960039,0.0769090454075157,1.97035358728859,1665673.81787078,3.3928973409089,27.901226461855,3.91324507714974,32.1802654568515,9.16761997255895,75.3892072968146
19637,0.0267482,0.0278142,0.9616742,0.9999933,1.57886185799754,0.0278141913003367,158.086068817119,158.087047080886,56.7646149028402,1.00980968677251,0.0175912391569409,0.988859514310117,1739697.71006039,1.63837023107481,58.904111695492,1.77234432104612,63.7208647164466,4.59142974744089,165.075076167514
\end{filecontents*}
\pgfplotstabletypeset[columns={dof,Qees,Qeex,theta,EffQoI,QeesSPR,thetaSPR,thetaQoISPR},
 columns/dof/.style
    ={column name=dof, column type=r, int detect},
 columns/theta/.style 
    ={column name=$\theta$,fixed,dec sep align, zerofill,precision=4},
 columns/EffQoI/.style 
    ={column name=$\theta_{QoI}$,fixed,dec sep align,zerofill,precision=5},
 columns/Qees/.style  
    ={column name={$\mathcal{E}$},fixed,dec sep align, zerofill,precision=4},
 columns/Qeex/.style    
    ={column name=$Q(\vm{e})$,fixed,dec sep align, zerofill,precision=4},
 columns/QeesSPR/.style  
    ={column name={$\mathcal{E}^{\dagger}$},fixed,dec sep align, zerofill,precision=4},
 columns/thetaSPR/.style 
    ={column name=$\theta^{\dagger}$,fixed,dec sep align, zerofill,precision=4},
 columns/thetaQoISPR/.style 
    ={column name=$\theta_{QoI}^{\dagger}$,fixed,dec sep align,zerofill,precision=5},
 col sep= comma]
{Data/CRACK_MI.csv} 
\label{tab:modeI}
\end{table}

\begin{table}
\centering
\caption{Stress intensity factor $K_{\rm II}$ as QoI. $^{\dagger}$Results using the standard SPR recovery.}
\begin{filecontents*}{Data/CRACK_MII.csv}
dof,Qees,Qeex,theta,EffQoI,QeesSPR,QeSPR,Q(u^h)SPR,F_d*U_pSPR,thetaSPR,thetaQoISPR,etaSPR,etaesSPR,K_dSPR,E2SPR,theta_2SPR,E3SPR,theta_3SPR,E4SPR,theta_4SPR
351,2.28859,1.12725,2.030241,1.007345,5.60519287091937,1.12725030203151,156.986632706387,157.006290211513,4.97244743320785,1.02832099549822,0.712935689506462,3.44740139980426,-17.8639523554907,6.55436854058929,5.81447485866902,13.3145634592416,11.811541265721,29.7210090865014,26.3659357934424
1289,0.282732,0.290166,0.974381,0.999953,2.22026826168592,0.290166110685647,157.823716897733,157.8363406537,7.65171458665365,1.0122070378279,0.183517161911833,1.3872862885002,0.419437924886507,3.19948895655133,11.0264046652144,5.26957974830266,18.1605623615071,15.3415523511434,52.8716200347867
4973,0.069965,0.068948,1.014748,1.000006,1.10783915917861,0.0689482876932175,158.044934720726,158.048214848245,16.0676819721454,1.0065705227885,0.0436067259758248,0.696085360104735,0.0663785689694306,1.37992122938894,20.0138578571936,1.94428371308624,28.199158791836,7.75808364522269,112.52032363359
19637,0.016755,0.015875,1.055441,1.000006,0.550150764155638,0.0158745266632536,158.098008481756,158.098825708199,34.6561995721818,1.00337905962036,0.0100399322065908,0.346774123803656,0.00422483900627953,0.622851784931754,39.2359279835116,0.766403577797591,48.2788302325677,3.88801533533471,244.921654535672
\end{filecontents*}
\pgfplotstabletypeset[columns={dof,Qees,Qeex,theta,EffQoI,QeesSPR,thetaSPR,thetaQoISPR},
 columns/dof/.style
    ={column name=dof, column type=r, int detect},
 columns/theta/.style 
    ={column name=$\theta$,fixed,dec sep align, zerofill,precision=4},
 columns/EffQoI/.style 
    ={column name=$\theta_{QoI}$,fixed,dec sep align,zerofill,precision=5},
 columns/Qees/.style  
    ={column name={$\mathcal{E}$},fixed,dec sep align, zerofill,precision=4},
 columns/Qeex/.style    
    ={column name=$Q(\vm{e})$,fixed,dec sep align, zerofill,precision=4},
 columns/QeesSPR/.style  
    ={column name={$\mathcal{E}^{\dagger}$},fixed,dec sep align, zerofill,precision=4},
 columns/thetaSPR/.style 
    ={column name=$\theta^{\dagger}$,fixed,dec sep align, zerofill,precision=4},
 columns/thetaQoISPR/.style 
    ={column name=$\theta_{QoI}^{\dagger}$,fixed,dec sep align,zerofill,precision=5},
 col sep= comma]
{Data/CRACK_MII.csv} 
\label{tab:modeII}
\end{table}

% \begin{table}[htb]
% \caption{Stress intensity factor $K_{\rm I}$ as QoI.}
% \centering
% \begin{tabular}{rcccc}
% \hline
% ndof & $Q(\vm{e}_{es})$ & $Q(\vm{e})$ & $\theta$ & $\theta_{QoI}$ \\ \hline
% 351 & 2.5264002 & 2.1144455 & 1.1948288 & 1.0026054 \\ 
% 1289 & 0.4821894 & 0.5146442 & 0.9369373 & 0.9997947 \\ 
% 4973 & 0.1140116 & 0.1216039 & 0.9375659 & 0.9999520 \\ 
% 19637 & 0.0267482 & 0.0278142 & 0.9616742 & 0.9999933 \\ \hline
% \end{tabular}
% \label{tab:modeI}
% \end{table}
% 
% \begin{table}[htb]
% \caption{Stress intensity factor $K_{\rm II}$ as QoI.}
% \centering
% \begin{tabular}{rcccc}
% \hline
% ndof & $Q(\vm{e}_{es})$ & $Q(\vm{e})$ & $\theta$ & $\theta_{QoI}$ \\ \hline
% 351 & 2.2885898 & 1.1272503 & 2.0302410 & 1.0073450 \\ 
% 1289 & 0.2827324 & 0.2901661 & 0.9743812 & 0.9999530 \\ 
% 4973 & 0.0699652 & 0.0689483 & 1.0147483 & 1.0000064 \\ 
% 19637 & 0.0167546 & 0.0158745 & 1.0554407 & 1.0000056 \\     \hline
% \end{tabular}
% \label{tab:modeII}
% \end{table}

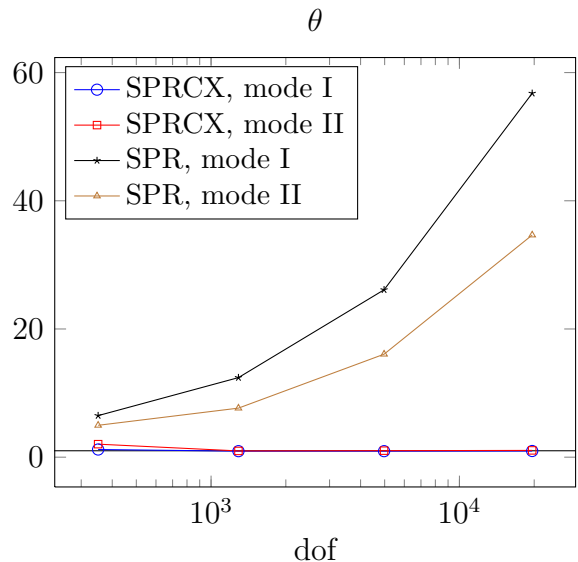
\begin{figure}[htb]
      \centering
 \begin{tikzpicture}%[scale=0.74]
    %\tikzset{/pgfplots/no markers}
    \begin{semilogxaxis}
    [
    %width=0.6\textwidth,
    title={$\theta$},
    xlabel={dof},
    legend style={cells={anchor=west}, font=\small},
    legend style={at={(0.02,0.98)},anchor=north west},
    cycle list name=ageplot,table/col sep=comma
    ]
    \draw (axis cs: 0,1)--(axis cs: 1e8,1); 
    \addplot table[x=dof,y=theta]{Data/CRACK_MI.csv};
     \addplot table[x=dof,y=theta]{Data/CRACK_MII.csv};
     \addplot table[x=dof,y=theta]{Data/SPRmodeI.csv};
     \addplot table[x=dof,y=theta]{Data/SPRmodeII.csv};
     \legend{
     {SPRCX, mode I},
     {SPRCX, mode II},
     {SPR, mode I},
     {SPR, mode II}}
    \end{semilogxaxis}
\end{tikzpicture}
    \caption{Evolution of the effectivity index $\theta$ for the SPRCX and SPR.
    Mode I and mode II loading conditions.}
\label{fig:thetaSPRCXvsSPR}
\end{figure}

Figure \ref{fig:thetaSPRCXvsSPR} compares the results of the proposed SPR-CX recovery with the standard SPR technique. In particular, the SPR cannot properly recover singular fields, thus, the error estimate provided by the technique does not converge to the exact error \cite{gonzalezrodenas2012}. This behaviour is similar for the two loading modes. 

% Figures \ref{fig:SPR-sigma_d}, \ref{fig:SPR-sigma_p} show 
% the $\sigma_{yy}$ stress component of the FE and recovered stress using the standard SPR technique for the primal and dual problems under mode I. In this case, the stress is smoothed at the singularity, as it is more clearly seen in Figure \ref{fig:SPR-sigma_d}. The SPR technique cannot describe the high gradients of the solution, which traduces in a lost of accuracy when compared with the proposed technique, especially in the vicinity of the crack tip where the contributions to the global error are considerable.
% 
% \begin{figure}[htb]
%     \centering
%     \includegraphics {Figures/SPR-2_sigma2h_p}
%     \includegraphics {Figures/SPR-2_sigma2s_p}    
%     \caption{FE (left) and recovered fields (right) $\sigma_{yy}$ for the primal problem using SPR.}
%     \label{fig:SPR-sigma_d}
% \end{figure}
% 
% \begin{figure}[htb]
%     \centering
%     \includegraphics {Figures/SPR-2_sigma2h_d}
%     \includegraphics {Figures/SPR-2_sigma2s_d}    
%     \caption{FE (left) and recovered fields (right) $\sigma_{yy}$ for the dual problem using SPR.}
%     \label{fig:SPR-sigma_p}
% \end{figure}

\begin{figure}[htb]
    \centering
    \includegraphics {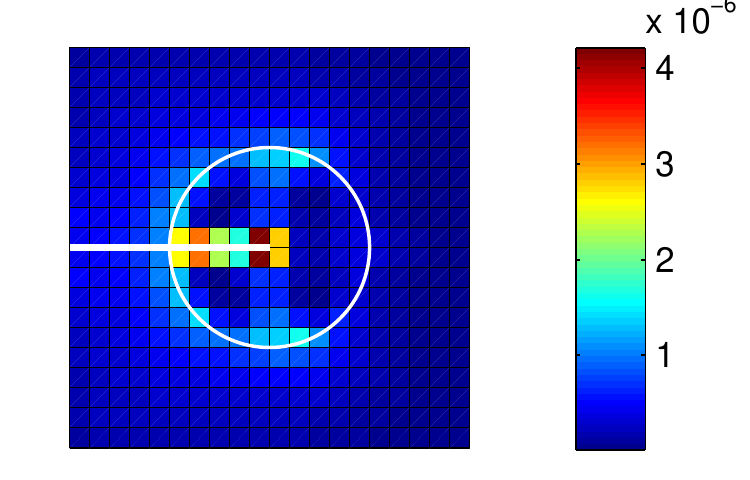}
    \includegraphics {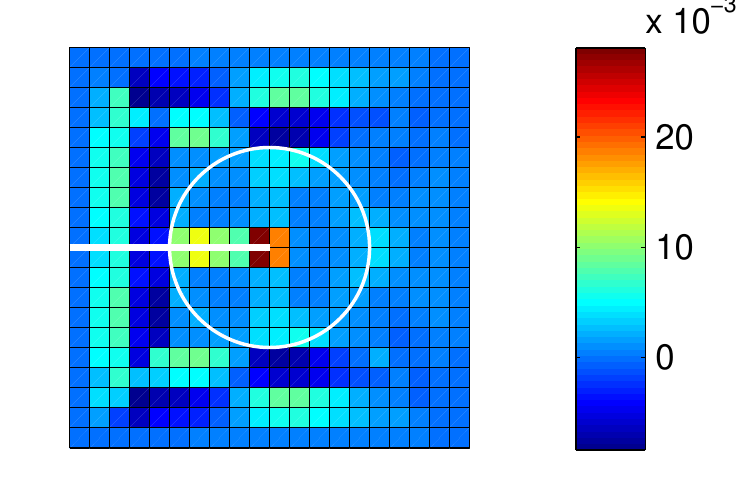}    
    \caption{Distribution of the estimated error considering the error in energy norm (top) and the error in the quantity of interest $\mathcal{E}$ (bottom).}
    \label{fig:ErrorDistribution}
\end{figure}

In Figure \ref{fig:ErrorDistribution} we represent the distribution of the estimated error for the second mesh of the sequence for the error in energy norm $\norm{\vm{e}_{es}}$ and the error considering the quantity of interest $\mathcal{E}$. This error distribution might guide the refinement during the adaptivity procedure. The approach based on energy norm estimates that the most critical part is located in the vicinity of the singular point whilst the goal oriented approach also considers the domain where the information of the QoI is extracted.

%TODO Include the Dual? why is so high?

%TODO include stresses in DoI

\section{Conclusions and future work.}
\label{sec:conclusions}

We have presented a locally equilibrated recovery procedure for goal oriented error estimation in XFEM. We have considered as the design parameter the generalised stress intensity factor that characterises the solution of singular problems in the context of linear elastic fracture mechanics.  The technique proposes the use of a stress recovery that provides locally equilibrated stress fields for both the primal and the dual problem. 

To formulate the dual problem we consider the linear equivalent domain integral representing $K$ to obtain the applied loads of the dual FE approximation. To perform the recovery of the primal and dual solutions we consider three main ideas: (i) enforcement of the internal equilibrium equation, (ii) enforcement of boundary equilibrium and (iii) splitting of the stress field into singular and smooth parts. 

The proposed technique has been tested with problems under different loading conditions. The obtained results show that the error estimator accurately captures the exact error in the evaluation of the stress intensity factor.

\section{Acknowledgements}

This work was supported by the EPSRC grant EP/G042705/1 ``Increased Reliability for Industrially Relevant Automatic Crack Growth Simulation with the eXtended Finite Element Method''. St\'ephane Bordas also thanks partial funding for his time provided by the European Research Council Starting Independent Research Grant (ERC Stg grant agreement No. 279578) ``RealTCut Towards real time multiscale simulation of cutting in non-linear materials with applications to surgical simulation and computer guided surgery''. This work has been carried out within the framework of the research project DPI2010-20542 of the Ministerio de Ciencia e Innovaci\'{o}n (Spain). The financial support of the FPU program (AP2008-01086), the funding from Universitat Polit\`{e}cnica de Val\`{e}ncia and Generalitat Valenciana (PROMETEO/2012/023) are also acknowledged. 

%\begin{thebibliography}{99}
%\end{thebibliography}
\bibliographystyle{wileyj-nourldoi}
\bibliography{library}

\begin{thebibliography}{10}
\providecommand{\url}[1]{\texttt{#1}}
\providecommand{\urlprefix}{URL }
\expandafter\ifx\csname urlstyle\endcsname\relax
  \providecommand{\doi}[1]{doi:\discretionary{}{}{}#1}\else
  \providecommand{\doi}{doi:\discretionary{}{}{}\begingroup
  \urlstyle{rm}\Url}\fi

\bibitem{ainsworthoden2000}
Ainsworth M, Oden JT. \emph{{A posteriori Error Estimation in Finite Element
  Analysis}}. John Wiley \& Sons: Chichester, 2000.

\bibitem{ladevezepelle2005}
Ladev\`{e}ze P, Pelle JP. \emph{{Mastering calculations in linear and nonlinear
  mechanics}}. Mechanical engineering series, Springer Science, 2005.

\bibitem{babuskarheinboldt1978}
Babu\v{s}ka I, Rheinboldt WC. {A-posteriori error estimates for the finite
  element method}. \emph{International Journal for Numerical Methods in
  Engineering}  1978; \textbf{12}(10):1597--1615.

\bibitem{zienkiewiczzhu1987}
Zienkiewicz OC, Zhu JZ. {A simple error estimator and adaptive procedure for
  practical engineering analysis}. \emph{International Journal for Numerical
  Methods in Engineering}  1987; \textbf{24}(2):337--357.

\bibitem{pereiraalmeida1999}
Pereira OJBA, de~Almeida JPM, Maunder EAW. {Adaptive methods for hybrid
  equilibrium finite element models}. \emph{Computer Methods in Applied
  Mechanics and Engineering}  Jul 1999; \textbf{176}(1-4):19--39.

\bibitem{paraschivoiuperaire1997}
Paraschivoiu M, Peraire J, Patera AT. {A posteriori finite element bounds for
  linear-functional outputs of elliptic partial differential equations}.
  \emph{Computer Methods in Applied Mechanics and Engineering}  1997;
  \textbf{150}(1-4):289--312.

\bibitem{ladevezerougeot1999}
Ladev\`{e}ze P, Rougeot P, Blanchard P, Moreau JP. {Local error estimators for
  finite element linear analysis}. \emph{Computer Methods in Applied Mechanics
  and Engineering}  1999; \textbf{176}(1-4):231--246.

\bibitem{odenprudhomme2001}
Oden JT, Prudhomme S. {Goal-oriented error estimation and adaptivity for the
  finite element method}. \emph{Computers \& Mathematics with Applications}
  Mar 2001; \textbf{41}(5-6):735--756.

\bibitem{cirakramm1998}
Cirak F, Ramm E. {A posteriori error estimation and adaptivity for linear
  elasticity using the reciprocal theorem}. \emph{Computer Methods in Applied
  Mechanics and Engineering}  Apr 1998; \textbf{156}(1-4):351--362.

\bibitem{ruterstein2006}
R\"{u}ter M, Stein E. {Goal-oriented a posteriori error estimates in linear
  elastic fracture mechanics}. \emph{Computer Methods in Applied Mechanics and
  Engineering}  Jan 2006; \textbf{195}(4-6):251--278.

\bibitem{almeidapereira2006}
de~Almeida JPM, Pereira OJBA. {Upper bounds of the error in local quantities
  using equilibrated and compatible finite element solutions for linear elastic
  problems}. \emph{Computer Methods in Applied Mechanics and Engineering}  Jan
  2006; \textbf{195}(4-6):279--296.

\bibitem{moesdolbow1999}
Mo\"{e}s N, Dolbow J, Belytschko T. {A finite element method for crack growth
  without remeshing}. \emph{International Journal for Numerical Methods in
  Engineering}  Sep 1999; \textbf{46}(1):131--150.

\bibitem{bordasduflot2007}
Bordas SPA, Duflot M. {Derivative recovery and a posteriori error estimate for
  extended finite elements}. \emph{Computer Methods in Applied Mechanics and
  Engineering}  2007; \textbf{196}(35-36):3381--3399.

\bibitem{duflotbordas2008}
Duflot M, Bordas SPA. {A posteriori error estimation for extended finite
  elements by an extended global recovery}. \emph{International Journal for
  Numerical Methods in Engineering}  2008; \textbf{76}:1123--1138.

\bibitem{rodenasgonzalez2008}
R\'{o}denas JJ, Gonz\'{a}lez-Estrada OA, Taranc\'{o}n JE, Fuenmayor FJ. {A
  recovery-type error estimator for the extended finite element method based on
  singular+smooth stress field splitting}. \emph{International Journal for
  Numerical Methods in Engineering}  2008; \textbf{76}(4):545--571.

\bibitem{rodenasgonzalez2010}
R\'{o}denas JJ, Gonz\'{a}lez-Estrada OA, D\'{\i}ez P, Fuenmayor FJ. {Accurate
  recovery-based upper error bounds for the extended finite element framework}.
  \emph{Computer Methods in Applied Mechanics and Engineering}  2010;
  \textbf{199}(37-40):2607--2621.

\bibitem{stroubouliszhang2006}
Strouboulis T, Zhang L, Wang D, Babu\v{s}ka I. {A posteriori error estimation
  for generalized finite element methods}. \emph{Computer Methods in Applied
  Mechanics and Engineering}  2006; \textbf{195}(9-12):852--879.

\bibitem{panetierladeveze2010}
Panetier J, Ladev\`{e}ze P, Chamoin L. {Strict and effective bounds in
  goal-oriented error estimation applied to fracture mechanics problems solved
  with XFEM}. \emph{International Journal for Numerical Methods in Engineering}
   2010; \textbf{81}(6):671--700.

\bibitem{gerasimovruter2012}
Gerasimov T, R\"{u}ter M, Stein E. {An explicit residual-type error estimator
  for Q 1 -quadrilateral extended finite element method in two-dimensional
  linear elastic fracture mechanics}. \emph{International Journal for Numerical
  Methods in Engineering}  2012; \textbf{90}(April):1118--1155.

\bibitem{panetierladeveze2009}
Panetier J, Ladev\`{e}ze P, Louf F. {Strict bounds for computed stress
  intensity factors}. \emph{Computers \& Structures}  Aug 2009;
  \textbf{87}(15-16):1015--1021.

\bibitem{rutergerasimov2013}
R\"{u}ter M, Gerasimov T, Stein E. {Goal-oriented explicit residual-type error
  estimates in XFEM}. \emph{Computational Mechanics}  2013;
  \textbf{52}(2):361--376.

\bibitem{pannachetsluys2009}
Pannachet T, Sluys LJ, Askes H. {Error estimation and adaptivity for
  discontinuous failure}. \emph{International Journal for Numerical Methods in
  Engineering}  2009; \textbf{78}(5):528--563.

\bibitem{gonzalezrodenas2012}
Gonz\'{a}lez-Estrada OA, R\'{o}denas JJ, Bordas SPA, Duflot M, Kerfriden P,
  Giner E. {On the role of enrichment and statical admissibility of recovered
  fields in a-posteriori error estimation for enriched finite element methods}.
  \emph{Engineering Computations}  2012; \textbf{29}(8).

\bibitem{diezrodenas2007}
D\'{\i}ez P, R\'{o}denas JJ, Zienkiewicz OC. {Equilibrated patch recovery error
  estimates: simple and accurate upper bounds of the error}.
  \emph{International Journal for Numerical Methods in Engineering}  2007;
  \textbf{69}(10):2075--2098.

\bibitem{rodenastur2007}
R\'{o}denas JJ, Tur M, Fuenmayor FJ, Vercher A. {Improvement of the
  superconvergent patch recovery technique by the use of constraint equations:
  the SPR-C technique}. \emph{International Journal for Numerical Methods in
  Engineering}  2007; \textbf{70}(6):705--727.

\bibitem{blackerbelytschko1994}
Blacker T, Belytschko T. {Superconvergent patch recovery with equilibrium and
  conjoint interpolant enhancements}. \emph{International Journal for Numerical
  Methods in Engineering}  1994; \textbf{37}(3):517--536.

\bibitem{szabobabuska1991}
Szab\'{o} BA, Babu\v{s}ka I. \emph{{Finite Element Analysis}}. John Wiley \&
  Sons: New York, 1991.

\bibitem{ginertur2009b}
Giner E, Tur M, Fuenmayor FJ. {A domain integral for the calculation of
  generalized stress intensity factors in sliding complete contacts}.
  \emph{International Journal of Solids and Structures}  Feb 2009;
  \textbf{46}(3-4):938--951.

\bibitem{zienkiewiczzhu1992}
Zienkiewicz OC, Zhu JZ. {The superconvergent patch recovery and a posteriori
  error estimates. Part 1: The recovery technique}. \emph{International Journal
  for Numerical Methods in Engineering}  1992; \textbf{33}(7):1331--1364.

\bibitem{rodenas2005}
R\'{o}denas JJ. {Goal Oriented Adaptivity: Una introducci\'{o}n a trav\'{e}s
  del problema el\'{a}stico lineal}. \emph{Technical {R}eport}, CIMNE, PI274,
  Barcelona, Spain 2005.

\bibitem{gonzalezrodenas2011b}
Gonz\'{a}lez-Estrada OA, R\'{o}denas JJ, Nadal E, Bordas SPA, Kerfriden P.
  {Equilibrated patch recovery for accurate evaluation of upper error bounds in
  quantities of interest}. \emph{Adaptive Modeling and Simulation. Proceedings
  of V ADMOS 2011}, Aubry D, D\'{\i}ez P, Tie B, Par\'{e}s N (eds.), CINME:
  Paris, 2011.

\bibitem{verdugodiez2011}
Verdugo F, D\'{\i}ez P, Casadei F. {Natural quantities of interest in linear
  elastodynamics for goal oriented error estimation and adaptivity}.
  \emph{Adaptive Modeling and Simulation. Proceedings of V ADMOS 2011}, Aubry
  D, D\'{\i}ez P, Tie B, Par\'{e}s N (eds.), CIMNE: Paris, 2011.

\bibitem{ginerfuenmayor2005}
Giner E, Fuenmayor FJ, Baeza L, Taranc\'{o}n JE. {Error estimation for the
  finite element evaluation of GI and GII in mixed-mode linear elastic fracture
  mechanics}. \emph{Finite Elements in Analysis and Design}  2005;
  \textbf{41}(11-12):1079--1104.

\bibitem{bechetminnebo2005}
B\'{e}chet E, Minnebo H, Mo\"{e}s N, Burgardt B. {Improved implementation and
  robustness study of the X-FEM for stress analysis around cracks}.
  \emph{International Journal for Numerical Methods in Engineering}  Oct 2005;
  \textbf{64}(8):1033--1056.

\bibitem{chessawang2003}
Chessa J, Wang H, Belytschko T. {On the construction of blending elements for
  local partition of unity enriched finite elements}. \emph{International
  Journal for Numerical Methods in Engineering}  Jun 2003;
  \textbf{57}(7):1015--1038.

\bibitem{graciewang2008}
Gracie R, Wang H, Belytschko T. {Blending in the extended finite element method
  by discontinuous Galerkin and assumed strain methods}. \emph{International
  Journal for Numerical Methods in Engineering}  2008;
  \textbf{74}(11):1645--1669.

\bibitem{fries2008}
Fries T. {A corrected XFEM approximation without problems in blending
  elements}. \emph{International Journal for Numerical Methods in Engineering}
  2008; \textbf{75}(5):503--532.

\bibitem{taranconvercher2009}
Taranc\'{o}n JE, Vercher A, Giner E, Fuenmayor FJ. {Enhanced blending elements
  for XFEM applied to linear elastic fracture mechanics}. \emph{International
  Journal for Numerical Methods in Engineering}  2009; \textbf{77}(1):126--148.

\end{thebibliography}

\end{document}